\documentclass[12pt,reqno]{amsart}

\usepackage{amsmath}
\usepackage{amsthm}
\usepackage{amssymb}
\usepackage{graphicx}
\usepackage{fullpage}

\newcommand{\N}{\mathbb{N}}
\newcommand{\R}{\mathbb{R}}
\renewcommand{\S}{\mathbb{S}}

\newcommand{\B}{\mathbb{B}}

\renewcommand{\d}{\mathrm{d}}
\newcommand{\mc}{\mathcal}

\newcommand{\I}{\mathrm{i}\,}
\newcommand{\wto}{\rightharpoonup}
\newcommand{\supp}{\operatorname{supp}}

\newtheorem{lemma}{Lemma}[section]
\newtheorem{proposition}[lemma]{Proposition}
\newtheorem{theorem}[lemma]{Theorem}
\newtheorem{corollary}[lemma]{Corollary}
\theoremstyle{remark}

\theoremstyle{definition}
\newtheorem{definition}[lemma]{Definition}

\numberwithin{equation}{section}
\numberwithin{table}{section}

\author{Pawe\l{} Biernat}
\address{Rheinische Friedrich-Wilhelms-Universit\"at Bonn, Mathematisches Institut, Endenicher Allee 60, D-53115 Bonn, Germany}
\email{pawel.biernat@gmail.com}

\author{Roland Donninger}
\address{Rheinische Friedrich-Wilhelms-Universit\"at Bonn,
Mathematisches Institut, Endenicher Allee 60, D-53115 Bonn, Germany}
\address{Universit\"at Wien, Fakult\"at f\"ur Mathematik, Oskar-Morgenstern-Platz 1, A-1090 Vienna, Austria}
\email{donninge@math.uni-bonn.de}

\author{Birgit Sch\"orkhuber}
\address{Universit\"at Wien, Fakult\"at f\"ur Mathematik, Oskar-Morgenstern-Platz 1, A-1090 Vienna, Austria}
\email{birgit.schoerkhuber@univie.ac.at}

\thanks{Roland Donninger is supported by the Alexander von Humboldt Foundation via
a Sofja Kovalevskaja Award endowed by the German Federal Ministry of Education
and Research. Birgit Sch\"orkhuber is supported by the Austrian Science Fund
(FWF) via the Hertha Firnberg Program, Project Nr. T 739-N25. Partial support by the Deutsche Forschungsgemeinschaft 
(DFG), CRC 1060 'The Mathematics of Emergent Effects', is also gratefully acknowledged.}

\title{Stable self-similar blowup in the supercritical heat flow of harmonic maps}

\begin{document}
\begin{abstract}
We consider the heat flow of corotational harmonic maps from $\R^3$ to the three-sphere and prove the nonlinear asymptotic stability of a particular self-similar shrinker that is not known in closed form. Our method provides a novel, systematic, robust, and constructive approach to the stability analysis of self-similar blowup in parabolic evolution equations. In particular, we completely avoid using delicate Lyapunov functionals, monotonicity formulas, indirect arguments, or fragile parabolic structure like the maximum principle. As a matter of fact, our approach reduces the nonlinear stability analysis of self-similar shrinkers to the spectral analysis of the associated self-adjoint linearized operators.
\end{abstract}

\maketitle

\section{Introduction}
\noindent Let $(M,g)$ and $(N, h)$ be Riemannian manifolds with metrics $g$ and $h$, respectively. A map $U: M\to N$ is called \emph{harmonic} if it is a critical point of the functional
\[ \mc S(U):=\int_M g^{jk}\partial_j U^a\partial_k U^b h_{ab}\circ U, \] 
where we employ Einstein's summation convention throughout.
Note that $\mc S(U)$ is a natural generalization of the Dirichlet energy. 
The Euler-Lagrange equations associated to $\mc S$ are
\[ \Delta_M U^a-g^{jk}\Gamma^a_{bc}(U)\partial_j U^b\partial_k U^c =0, \]
where $\Gamma^a_{bc}$ are the Christoffel symbols on the target manifold $N$
and 
\[ \Delta_M=\frac{1}{\sqrt{\det g}}\partial_j \left (\sqrt{\det g}\,g^{jk}\partial_k\right ) \] is the Laplace-Beltrami operator on $M$.
The study of harmonic maps is a classical subject in geometric analysis, see e.g.~\cite{EelSam64, Smi75, SacUhl81, Str85, Str92, EelRat93, EelLem95, SchYau97, Hel02, LinCha08}.
The basic mathematical questions concern the existence and, ideally, the classification of harmonic maps.
A standard tool in this respect is the associated \emph{heat flow}, i.e., one considers a one-parameter family $\{U_t: t\geq 0\}$ of maps from $M$ to $N$ that evolve according to the heat equation
\[ \partial_t U_t^a=\Delta_M U_t^a-g^{jk}\Gamma^a_{bc}(U_t)\partial_j U_t^b\partial_k U_t^c. \] 
The idea then is to take an arbitrary map $U_0: M\to N$ as initial data at $t=0$ and due to the regularizing effects of the heat flow, the solution $U_t$ is expected to converge to an equilibrium as $t\to\infty$.
In other words, the heat flow is supposed to deform arbitrary maps into harmonic ones. Indeed, this strategy works well under certain curvature assumptions as is demonstrated in the classical paper \cite{EelSam64}.
In the general case, however, the flow tends to form singularities (or ``blow up'') 
in finite time \cite{CorGhi89, CheDin90, ChaDinYe92, Qin95, QinTia97, Fan99, Gas02, BerBouHulKin03, Top04, GuaGusTsa09, BieBiz11, RapSch13, RapSch14, Bie15}. This is a severe obstruction which can only be overcome if one is able to continue the flow past the singularity in a well-defined manner. 
Such a construction is a challenging endeavor which presupposes a detailed understanding of possible blowup scenarios. 
Naturally, one is mainly interested in blowup behavior that is stable under small perturbations of the initial data.

In this paper we are interested in singularity formation in the heat flow of harmonic maps $U: \S^d\to\S^d$. As it turns out, the blowup is a local phenomenon and the curvature of the base manifold is irrelevant for the asymptotic behavior near the singularity. Consequently, we may equally well consider maps $U: \R^d\to \S^d$, cf.~\cite{Str88, Gas02}. Furthermore, we restrict ourselves to the case $d=3$ and assume corotational symmetry. That is to say, we choose standard spherical coordinates $(r,\theta,\varphi)$ on $\R^3$, hyperspherical coordinates on $\S^3$, and make the ansatz $U(r,\theta,\varphi)=(u(r),\theta,\varphi)$ for the map $U: \R^3\to\S^3$.
Under this symmetry reduction, the Euler-Lagrange equations associated to the functional $\mc S$ reduce to a single nonlinear ordinary differential equation for $u$ which reads
\[ u''(r)+\frac{2}{r}u'(r)-\frac{\sin(2u(r))}{r^2}=0,\qquad r\geq 0. \]
In order to obtain the associated heat flow, we introduce an artificial time dependence and consider the Cauchy problem for the equation
\begin{equation}
\label{eq:hmhf} \partial_t u(r,t)-\partial_r^2 u(r,t)-\frac{2}{r}\partial_r u(r,t)+\frac{\sin(2u(r,t))}{r^2}=0. 
\end{equation}
Our main result shows the existence of a stable self-similar blowup scenario for Eq.~\eqref{eq:hmhf}. 
For the precise formulation we introduce the following function space.

\begin{definition}
Let
\[ \tilde Y:=\{h\in C^\infty_c([0,\infty)): h^{(2k)}(0)=0 \mbox{ for all }k\in \N_0\} \]
and set
\[ \|h\|_Y:=\||\cdot|^{-1}h(|\cdot|)\|_{\dot H^2(\R^5)}+\||\cdot|^{-1}h(|\cdot|)\|_{\dot H^4(\R^5)}. \]
The Banach space $Y$ is defined as the completion of $\tilde Y$ with respect to $\|\cdot\|_Y$.
\end{definition}

\begin{theorem}
\label{thm:main}
There exists an $f_0 \in C^\infty([0,\infty))\cap Y$ with $f_0>0$ on $(0,\infty)$ such that, for any $T_0>0$ and $t\in [0,T_0)$,
\[ u^*_{T_0}(r,t):=f_0\left (\frac{r}{\sqrt{T_0-t}} \right ) \]
is a solution to Eq.~\eqref{eq:hmhf}.
Furthermore, there exist $\delta,M,\omega_0>0$ such that the following holds.
For any $h\in Y$ satisfying $\|h\|_Y \leq \frac{\delta}{M^2}$, there exists a $T_h \in [T_0-\frac{\delta}{M},T_0+\frac{\delta}{M}]$ such that Eq.~\eqref{eq:hmhf} with initial data $u(r,0)=u^*_{T_0}(r,0)+h(r)$ has a unique solution $u_h$ that blows up at $t=T_h$ and converges to $u^*_{T_h}$ in the sense that
\[ \frac{\|u_h(\cdot,t)-u^*_{T_h}(\cdot,t)\|_Y}{\|u_{T_h}^*(\cdot, t)\|_Y}
\leq \delta(T_h-t)^{\omega_0} \]
for all $t\in [0,T_h)$.
In particular, the class $\{u^*_{T_0}: T_0>0\}$ of self-similar solutions is nonlinearly asymptotically stable under small perturbations of the initial data.
\end{theorem}

Some remarks are in order.
\begin{itemize}

\item The map $U: \R^3\to\S^3$ has values on the sphere and thus, there is no blowup in $L^\infty$. However, the self-similar solution $u_{T_0}^*$ blows up in $Y$. Indeed, a simple scaling argument shows
\[ \|u_{T_0}^*(\cdot,t)\|_Y\simeq (T_0-t)^{-\frac54} \]
for $t\in [0,T_0)$.

\item The blowup profile $f_0$ is constructed in the companion paper \cite{BieDon16} by a novel computer-assisted (but rigorous) method. It is not known in closed form. Furthermore, $f_0$ is not the only self-similar profile. In fact, there exist infinitely many self-similar solutions to Eq.~\eqref{eq:hmhf}, see \cite{Fan99}. To the knowledge of the authors, Theorem \ref{thm:main} is the first result on stable blowup with a nonunique blowup profile that is not known explicitly. 

\item The norm $\|\cdot\|_Y$ might look odd at first glance since it is based on homogeneous Sobolev spaces on $\R^5$ whereas Eq.~\eqref{eq:hmhf} is posed on $\R^3$.
However, if one sets $u(r,t)=rv(r,t)$, Eq.~\eqref{eq:hmhf} transforms into a radial heat equation on $\R^5$ for the function $v$. In addition, this transformation regularizes the nonlinearity at the center, see below. In this sense, the effective dimension of the problem is $5$ and it is natural to work with radial functions on $\R^5$.

\item In the formulation of Theorem \ref{thm:main} we do not specify the precise solution concept we are using. We will study Eq.~\eqref{eq:hmhf} in similarity coordinates by semigroup theory which yields a canonical notion of strong solution (which is actually called ``mild solution'' in semigroup theory). Since Eq.~\eqref{eq:hmhf} is parabolic, smoothing effects will kick in immediately and turn strong solutions into classical ones.

\item For obvious reasons, self-similar solutions of the form $f(\frac{r}{\sqrt{T_0-t}})$ are called \emph{shrinkers}. Since Eq.~\eqref{eq:hmhf} is not time-reversible, there is another, independent class of self-similar solutions, so-called \emph{expanders}, which take the form $f(\frac{r}{\sqrt{t-T_0}})$. The latter have also attracted considerable interest, in particular in connection with the question of unique continuation beyond blowup \cite{BieBiz11, GerRup11, GerGhoMiu16}, but they play no role in the present paper.
\end{itemize}

\subsection{Related results}
The analysis of harmonic maps is a vast subject that is impossible to review in this paper. We restrict ourselves to a brief discussion of recent blowup results that are directly related to our work and refer the reader to the monographs and survey articles \cite{Str92, EelRat93, EelLem95, SchYau97, Hel02, LinCha08} for the general background.
 
As already indicated, self-similar solutions for the corotational heat flow of harmonic maps $U:\R^d\to\S^d$ for $d\in \{3,4,5,6\}$ are constructed in \cite{Fan99, Gas02}.
Expanding self-similar solutions are studied in \cite{GerRup11}.
For $d\geq 7$, there are no self-similar shrinkers \cite{BizWas15} and the blowup is of a more complicated nature \cite{Bie15, BieSek16}. 
The case $d=2$ is of special interest since it is energy-critical and blowup takes place via shrinking of a soliton \cite{BerBouHulKin03, RapSch13, RapSch14}.
The unique continuation beyond blowup is investigated in \cite{BieBiz11, GerGhoMiu16}.
Needless to say, there are similar results for closely related problems like the Yang-Mills heat flow or the nonlinear heat equation, see the discussion in \cite{DonSch16} for a brief overview.
Of particular interest in this context is the recent paper \cite{ColRapSze16} which also considers self-similar blowup for a nonlinear heat equation with a blowup profile that is not known in closed form. In contrast to our result, however, the blowup studied in \cite{ColRapSze16} is highly unstable and the necessary spectral properties can be obtained by a perturbative argument.

\subsection{Outline of the proof}
The proof of Theorem \ref{thm:main} proceeds by a perturbative construction around the blowup solution $u_{T_0}^*$. We would like to emphasize that this is a robust approach that uses no structure other than the spectral stability of the self-similar profile $f_0$ which is established in \cite{BieDon16}. As a consequence, our method provides a universal framework for studying self-similar blowup in general parabolic evolution equations. We briefly outline the main steps.

\begin{itemize}
\item We consider Eq.~\eqref{eq:hmhf} with initial data $u(r,0)=u_{T_0}^*(r,0)+h(r)$. By time translation invariance we may assume $T_0=1$ and we introduce similarity coordinates $s=-\log(T-t)+\log T$, $y=\frac{r}{\sqrt{T-t}}$ which go back to \cite{GigKoh85}.
Here, $T>0$ is a free parameter which will be adjusted later.
Then we rescale the dependent variable $u$ in a suitable manner to obtain the evolution equation
\[ \partial_s \tilde w-\partial_y^2 \tilde w-\tfrac{4}{y}\partial_y \tilde w+\tfrac{y}{2}\partial_y \tilde w-\tfrac{2}{y^2}\tilde w+\tfrac12 \tilde w+\tfrac{\sin(y\tilde w)}{y^3}=0, \]
where $\tilde w=\tilde w(y,s)$, with initial data $\tilde w(y,0)=f_0(\sqrt Ty)/y+h(\sqrt T y)/y$.
This equation has the static solution $\tilde w(y,s)=f_0(y)/y$.
To study its stability, we make the ansatz $\tilde w(y,s)=f_0(y)/y+w(y,s)$ which leads to an evolution equation of the form
\begin{equation}
\label{eq:introw}
\left \{\begin{array}{l}
 \partial_s w(\cdot,s)=\hat{\mc L}w(\cdot,s)+\mc N(w(\cdot,s)) \\
 w(y,0)=f_0(\sqrt T y)/y-f_0(y)/y+h(\sqrt Ty)/y
 \end{array} \right .
 \end{equation}
for the perturbation $w$. The linear operator $\hat{\mc L}$ is given by
\[ \hat{\mc L}=\partial_y^2 +\tfrac{4}{y}\partial_y -\tfrac{y}{2}\partial_y-\tfrac12 -V_0(y) \]
with the potential $V_0(y)=\tfrac{2\cos(2f_0(y))-2}{y^2}$ and $\mc N$ denotes the nonlinear remainder.
In the spirit of standard local well-posedness theory we now try to solve Eq.~\eqref{eq:introw} by treating the nonlinear terms perturbatively. Consequently, we first have to understand the linearized equation that arises from \eqref{eq:introw} by dropping the nonlinear terms.

\item The operator $\hat{\mc L}$, interpreted as an operator acting on radial functions on $\R^5$, has a self-adjoint extension $\mc L$ on $L^2_\sigma(\R^5)$ with the weight $\sigma(x)=e^{-|x|^2/4}$.
Here we encounter the fundamental problem in studying self-similar blowup for parabolic equations:
In order to apply self-adjoint spectral theory, it seems necessary to study the evolution in Sobolev spaces with exponentially decaying weights. This, however, is impossible since one cannot control nonlinear terms in such spaces.

There are (at least) two ways around this issue. First, one can study the evolution in unweighted Sobolev spaces and rely on \emph{non}self-adjoint spectral theory.
This approach was chosen in \cite{DonSch16} for the study of the Yang-Mills heat flow.
In this paper we follow a different strategy which is based on the simple observation that in a certain sense the problem splits into a self-adjoint part on a compact domain, where the exponentially decaying weight is irrelevant, and a nonself-adjoint part on an unbounded domain which, however, is easy since the potential term is negligible there.
We remark that this is not a new discovery but a well-known phenomenon in parabolic problems, see e.g.~\cite{BriKup94, MerZaa97, TayZaa15, ColRapSze16}.
Somewhat paradoxically, we can therefore study the linearized evolution on \emph{unweighted} spaces by using \emph{self-adjoint} spectral theory in a \emph{weighted} space.

More precisely, we consider the semigroup $e^{s\mc L}$ on $L^2_\sigma(\R^5)$ generated by the self-adjoint operator $\mc L$. From \cite{BieDon16} we know that $\mc L$ has precisely one nonnegative eigenvalue $\lambda=1$ with eigenfunction $\psi_1$.
As usual, this instability is related to the freedom in choosing the parameter $T$ in the similarity coordinates.
From self-adjoint spectral theory we obtain the weighted decay estimate 
\[ \|e^{s\mc L}f\|_{L^2_\sigma}\lesssim e^{-c_0 s}\|f\|_{L^2_\sigma} \]for some constant $c_0>0$, provided $f\perp \psi_1$.
Similar bounds hold for higher Sobolev spaces with weights.
As a matter of fact, also on \emph{unweighted} homogeneous Sobolev spaces of sufficiently high degree we have decay, but a priori only for the \emph{free operator} $\mc L_0=\mc L-V_0$.
Indeed, an integration by parts shows
\[ (\Delta \mc L_0 f | \Delta f)_{L^2}\leq -\tfrac14 \|f\|_{L^2}^2 \]
on the unweighted $L^2$. Similar bounds hold for higher derivatives.
Consequently, by combining the unweighted bounds, the weighted decay, and the smallness of $V_0(y)$ for large $y$, we derive the unweighted decay 
\[ \|e^{s\mc L}f\|_X\lesssim e^{-\omega_0 s}\|f\|_X,\qquad f\perp \psi_1 \] for some $\omega_0>0$, where $X=\dot H^2(\R^5)\cap \dot H^4(\R^5)$.

\item From now on we follow the argument introduced in our earlier works \cite{Don11, DonSch12, DonSch14, Don14, DonSch15, DonSch16a} on self-similar blowup for wave-type equations. We first show that the nonlinearity is locally Lipschitz on $X$. This is not hard but requires at least some work due to the removable singularity of the nonlinearity at the center.
Then we employ Duhamel's principle to rewrite Eq.~\eqref{eq:introw} as
\begin{equation}
\label{eq:introphi} \phi(s)=e^{s\mc L}\mc U(h,T)+\int_0^s e^{(s-s')\mc L}\mc N(\phi(s'))\d s' 
\end{equation}
where $\phi(s)(y)=w(y,s)$ and $\mc U(h,T)$ is an abbreviation for the initial data.
In general, Eq.~\eqref{eq:introphi} does not have a global solution due to the unstable eigenvalue $1\in \sigma(\mc L)$. 
We deal with this issue by employing the Lyapunov-Perron method.
That is to say, we first suppress the instability by subtracting a correction term and instead of Eq.~\eqref{eq:introphi}, we consider the modified equation
\begin{equation}
\label{eq:intromod} \phi(s)=e^{s\mc L}\mc U(h,T)+\int_0^s e^{(s-s')\mc L}\mc N(\phi(s'))\d s' -e^s\mc C(\phi,\mc U(h,T)) 
\end{equation}
with
\[ \mc C(\phi,\mc U(h,T))=\mc P \mc U(h,T)+\int_0^\infty e^{-s'}\mc P\mc N(\phi(s'))\d s'. \]
Here, $\mc P$ is the orthogonal projection on the unstable subspace $\langle \psi_1\rangle$.
By a fixed point argument we show that for any small $h$ and $T$ close to $1$, Eq.~\eqref{eq:intromod} has a global (in $s$) solution $\phi_{h,T}$ that decays like the stable linear flow, i.e., $\|\phi_{h,T}(s)\|\lesssim e^{-\omega_0 s}$. 
In the final step we prove that for any small $h$, there exists a $T_h$ close to $1$ which makes the correction term vanish. In other words, $\phi_{h,T_h}$ is a solution to the original equation \eqref{eq:introphi}.

\end{itemize}

\section{Preliminary transformations}
\label{sec:pre}

\noindent The basic evolution equation is
\begin{equation}
\label{eq:u} 
\partial_t u(r,t)-\partial_r^2 u(r,t)-\frac{2}{r}\partial_r u(r,t)+\frac{\sin(2u(r,t))}{r^2}=0 
\end{equation}
where $r\geq 0$.
For any $T_0>0$, we have the self-similar solution 
\[ u_{T_0}^*(r,t)=f_0\left (\frac{r}{\sqrt{T_0-t}}\right ) \]
with $f_0$ constructed in \cite{BieDon16}.
Our goal is to study the evolution of small initial perturbations of $u_{T_0}^*$.
By time translation invariance, we may restrict ourselves to $T_0=1$.
Consequently, we consider the Cauchy problem
\begin{equation}
\left \{ \begin{array}{l}
\partial_t u(r,t)-\partial_r^2 u(r,t)-\frac{2}{r}\partial_r u(r,t)+\frac{\sin(2u(r,t))}{r^2}=0 \\
u(r,0)=u_1^*(r,0)+h(r)=f_0(r)+h(r)
\end{array} \right. 
\end{equation}
where $h$ is a free function.
In order to regularize the nonlinearity, it is useful to change variables according to $u(r,t)=rv(r,t)$. This yields
\begin{equation}
\label{eq:v}
\left \{ \begin{array}{l}
\partial_t v(r,t)-\partial_r^2 v(r,t)-\frac{4}{r}\partial_r v(r,t)-\frac{2}{r^2}v(r,t)+\frac{\sin(2rv(r,t))}{r^3}=0 \\
v(r,0)=f_0(r)/r+h(r)/r.
\end{array} \right. 
\end{equation}
Accordingly, we write $u_{T_0}^*(r,t)=rv_{T_0}^*(r,t)$ for the self-similar solution.
Now we switch to similarity coordinates
$s=-\log(T-t)+\log T$, $y=\frac{r}{\sqrt{T-t}}$ and define the new dependent variable $\tilde w$ by
\[ \tilde w(y,s):=\sqrt T e^{-s/2}v\left (\sqrt T y e^{-s/2}, T(1-e^{-s})\right ), \]
or, equivalently,
\[ v(r,t)=\frac{1}{\sqrt{T-t}}\tilde w\left (\frac{r}{\sqrt{T-t}}, -\log(T-t)+\log T\right ). \]
Here, $T>0$ is a free parameter that will be needed to account for the time translation invariance of the problem which introduces an artificial instability.
Eq.~\eqref{eq:v} transforms into
\begin{equation}
\label{eq:w}
\left \{ \begin{array}{l}
\partial_s \tilde w(y,s)-\partial_y^2 \tilde w(y,s)-\frac{4}{y}\partial_y \tilde w(y,s)+\frac{y}{2}\partial_y \tilde w(y,s)-\frac{2}{y^2}\tilde w(y,s)+\frac12 \tilde w(y,s)+\frac{\sin(2y\tilde w(y,s))}{y^3}=0 \\
\tilde w(y,0)=f_0(\sqrt T y)/y+h(\sqrt T y)/y.
\end{array} \right.
\end{equation}
Observe that the only trace of the parameter $T$ is in the initial data.
Furthermore, by construction, 
\begin{align*}
 \tilde w_T^*(y,s):&=\sqrt T e^{-s/2}v_T^*\left (\sqrt Tye^{-s/2},T(1-e^{-s})\right )=
\frac{1}{y}u_T^*\left (\sqrt Tye^{-s/2},T(1-e^{-s})\right) \\
&=f_0(y)/y 
\end{align*} is a static solution to Eq.~\eqref{eq:w}.
By making the ansatz $\tilde w(y,s)=f_0(y)/y+w(y,s)$, we rewrite Eq.~\eqref{eq:w} as 
\begin{equation}
\label{eq:full}
\left \{ \begin{array}{l}
 \partial_s w(y,s)=\hat{\mc L} w(y,s)+\hat{\mc N}(w(y,s)) \\
 w(y,0)=f_0(\sqrt T y)/y-f_0(y)/y+h(\sqrt Ty)/y
 \end{array} \right.
 \end{equation}
with the linear operator $\hat{\mc L}$ defined by
\begin{align}
\label{def:hatL}
\begin{split}
\hat {\mc L}&=\partial_y^2 +\frac{4}{y}\partial_y -\frac{y}{2}\partial_y-\frac12 -\frac{2\cos(2f_0(y))-2}{y^2} 
\end{split}
\end{align}
and the nonlinearity
\begin{equation}
\hat{\mc N}(w(y,s))=-\frac{1}{y^3}\left [\sin(2f_0(y)+2yw(y,s))-\sin(2f_0(y))-2y\cos(2f_0(y))w(y,s) \right ].
\end{equation}

\section{The linearized evolution}

\noindent In this section we study the linearized equation, i.e., we drop the nonlinearity in Eq.~\eqref{eq:full} and focus on 
\begin{equation}
\label{eq:lin}
 \partial_s w(y,s)=\hat{\mc L}w(y,s). 
 \end{equation}
Furthermore, we do not specify the initial data explicitly because their specific form is irrelevant for the linear theory.

Note that the operator $\hat{\mc L}$ contains the $5$-dimensional radial Laplacian and 
for the rest of this paper we actually find it convenient to switch to $5$-dimensional notation.
To this end, we define the operator
\[ \Lambda f(x):=\tfrac12 x\nabla f(x)+\tfrac12 f(x) \]
acting on functions $f: \R^5\to \R$.
In the following, the variable $x$ is used to denote an element of $\R^5$.
In this spirit we define the potential $V_0: \R^5\to\R$ by
\[ V_0(x):=-\frac{2\cos(2f_0(|x|))-2}{|x|^2}. \]
By \cite{BieDon16}, $f_0$ is odd\footnote{By this we mean that $f_0$ can be extended to all of $\R$ as a smooth, odd function. In other words, $f_0^{(2k)}(0)=0$ for all $k\in \N_0$.} and thus, $V_0\in C^\infty(\R^5)$, see \cite{Whi43}.
Now we define a differential operator $\tilde{\mc L}$ by
\[ \tilde{\mc L} f:=\Delta f-\Lambda f+V_0 f \]
where throughout, $\Delta$ denotes the Laplacian on $\R^5$.
Then we have
\[ \tilde{\mc L}f(x)=\tilde f''(|x|)+\frac{4}{|x|}\tilde f'(|x|)
-\frac{|x|}{2}\tilde f'(|x|)-\frac12 \tilde f(|x|)+V_0(x)\tilde f(|x|) \]
for all radial functions $f:\R^5\to\R$ with $f(x)=\tilde f(|x|)$.
Consequently, the linearized equation \eqref{eq:lin} can be written as
\begin{equation}
\label{eq:linop}
 \partial_s \phi(s)=\tilde{\mc L}\phi(s) 
 \end{equation}
where $\phi(s)(x)=w(|x|,s)$.
Formally, the solution of Eq.~\eqref{eq:linop} is given by $\phi(s)=e^{s\tilde{\mc L}}\phi(0)$. 
In the following, we make this rigorous.

\subsection{Basic semigroup theory}
As usual, for $\Omega\subset \R^d$ open and $w: \Omega\to [0,\infty)$ a weight function, we write
\[ (f|g)_{L^2_w(\Omega)}:=\int_\Omega f(x)g(x)w(x)\d x,\qquad \|f\|_{L^2_w(\Omega)}:=\sqrt{(f|f)_{L^2_w(\Omega)}} \]
and denote by $L^2_w(\Omega)$ the completion of $C^\infty_c(\Omega)$ with respect to $\|\cdot\|_{L^2_w(\Omega)}$.

We promote $\tilde{\mc L}$ to an unbounded linear operator on the Hilbert space 
\[ H:=\{f\in L^2_\sigma(\R^5): f \mbox{ radial}\} \] 
with $\sigma(x)=e^{-|x|^2/4}$, by specifying the domain $\mc D(\tilde{\mc L}):=\{f\in C^\infty_c(\R^5): f\mbox{ radial}\}$.

\begin{proposition}
\label{prop:gen}
The operator $\tilde{\mc L}: \mc D(\tilde{\mc L})\subset H\to H$ is essentially self-adjoint and the spectrum of its closure $\mc L$ satisfies $\sigma(\mc L)\cap [0,\infty)=\{1\}$.
The spectral point $1$ is a simple eigenvalue
and $\mc L$ generates a strongly continuous one-parameter semigroup $e^{s\mc L}$ on $H$. 
The function $\psi_1(x):=f_0'(|x|)/\|f_0'(|\cdot|)\|_{L^2_\sigma(\R^5)}$ is an eigenfunction of $\mc L$ with eigenvalue $1$. Moreover, there exists a constant $c_0>0$ such that
\[ \|e^{s\mc L}f\|_{L^2_\sigma(\R^5)}\leq e^{-c_0 s}\|f\|_{L^2_\sigma(\R^5)} \]
for all $f\in H$ satisfying $(f|\psi_1)_{L^2_\sigma(\R^5)}=0$
and all $s\geq 0$.
\end{proposition}

\begin{proof}
Via $\tilde f\mapsto |\S^4|^{-1/2}\tilde f(|\cdot|): L^2_\rho(0,\infty)\to H$ with the weight $\rho(y)=y^4 e^{-y^2/4}$, $\tilde{\mc L}$ is unitarily equivalent to the Sturm-Liouville operator
\begin{equation*}
\mc T \tilde f(y):=\frac{1}{\rho(y)}\frac{\d}{\d y}\left [\rho(y)\tilde f'(y)\right ]-\frac12 \tilde f(y)-\frac{2\cos(2f_0(y))-2}{y^2} \tilde f(y)
\end{equation*}
with domain $\mc D(\mc T):=\{\tilde f\in C^\infty_c([0,\infty)): \tilde f^{(2k)}(0)=0 \mbox{ for all }k\in \N_0\}$.
The equation
\begin{equation}
\label{eq:SL} 
\frac{1}{\rho(y)}\frac{\d}{\d y}\left [\rho(y)\tilde f'(y)\right ]=0 
\end{equation}
has the explicit solution
\[ \tilde f_1(y)=\int_1^y \rho(s)^{-1}\d s = \int_1^y s^{-4}e^{s^2/4}\d s .\]
For $y\in (0,1]$ we have
\[ |\tilde f_1(y)|=\int_y^1 s^{-4}e^{s^2/4}\d s\geq \int_y^1 s^{-4}\d s=\tfrac13 y^{-3}-\tfrac13 \]
and thus, $\tilde f_1\notin L^2_\rho(0,1)$.
Similarly, for $y\geq 1$,
\[ \tilde f_1(y)=2\int_1^y s^{-5}\partial_s e^{s^2/4}\d s\geq 2y^{-5}(e^{y^2/4}-e^{1/4}) \]
which implies $\tilde f_1\notin L^2_\rho(1,\infty)$.
By the Weyl alternative, the Sturm-Liouville operator defined by \eqref{eq:SL} is in the limit-point case at both endpoints and the Kato-Rellich theorem implies that $\mc T$ (and hence $\tilde{\mc L}$) is essentially self-adjoint, see e.g.~\cite{Tes14}.  

In fact, by $\tilde f\mapsto |\S^4|^{-1/2}\tilde f(|\cdot|)/|\cdot|: L^2_{\tilde \rho}(0,\infty)\to H$ with $\tilde \rho(y)=y^2 e^{-y^2/4}$, $\mc L$ is unitarily equivalent to the operator $-\mc A_0$ studied in \cite{BieDon16}. Consequently, from \cite{BieDon16} we obtain $\sigma(\mc L)\cap [0,\infty)=\{1\}$ with $1$ a simple eigenvalue.
The corresponding normalized eigenfunction is given by $\psi_1(x)=f_0'(|x|)/\|f_0(|\cdot|)\|_{L^2_\sigma(\R^5)}$, see \cite{BieDon16}.
Since $0\notin \sigma(\mc L)$, we obtain $-c_0:=\sup\sigma(\mc L)\setminus\{1\}<0$ and the self-adjointness of $\mc L$ implies the bound
\begin{equation}
\label{eq:diss}
 (\mc L f | f)_{L^2_\sigma(\R^5)}\leq -c_0 \|f\|_{L^2_\sigma(\R^5)}^2 
 \end{equation}
for all $f\in \mc D(\mc L)$ with $(f|\psi_1)_{L^2_\sigma(\R^5)}=0$.
From this, the stated bound on the semigroup $e^{s\mc L}$ follows.
\end{proof}

\subsection{Estimates in local Sobolev norms}

We upgrade the $L^2_\sigma$ bound on $e^{s\mc L}$ to a local $H^4$ bound.
In the following we use
\[ \|f\|_{G(\mc L)}:=\|\mc Lf\|_{L^2_\sigma(\R^5)}+\|f\|_{L^2_\sigma(\R^5)} \]
for $f\in \mc D(\mc L)$ to denote the graph norm of $\mc L$.
Furthermore, the letter $C$ (possibly with subscripts to indicate dependencies) denotes a positive constant that might change its value at each occurrence and $c_0>0$ is the constant from Proposition \ref{prop:gen}.
Finally, for $R>0$ we set
\[ \B^5_R:=\{x\in \R^5: |x|<R\}. \]

\begin{lemma}
\label{lem:H2B}
Let $f\in \mc D(\mc L)$ and $R\geq 1$.
Then $\nabla f, \Delta f\in L^2(\B^5_R)$ and 
we have the bound 
\[ \|\Delta f\|_{L^2(\B_R^5)}+\|\nabla f\|_{L^2(\B^5_R)}\leq C_R \|f\|_{G(\mc L)} \]
for all $R\geq 1$ and all $f\in \mc D(\mc L)$.
\end{lemma}

\begin{proof}
Let $f\in C^\infty_c(\R^5)$ and $R\geq 1$.
An integration by parts yields
\[ (\mc L f | f)_{L^2_\sigma(\R^5)}\leq -\|\nabla f\|_{L^2_\sigma(\R^5)}^2+C\|f\|_{L^2_\sigma(\R^5)}^2 \]
and we infer
\begin{equation}
\label{eq:nabf} \|\nabla f\|_{L^2(\B^5_R)}\leq C_R \|\nabla f\|_{L^2_\sigma(\R^5)}\leq C_R \|f\|_{G(\mc L)}. 
\end{equation}
Now let $f\in \mc D(\mc L)$. Since $C^\infty_c(\R^5)$ is a core for $\mc L$, there exists a sequence $(f_n)\subset C^\infty_c(\R^5)$ such that $f_n \to f$ in the graph norm $\|\cdot\|_{G(\mc L)}$. Consequently, Eq.~\eqref{eq:nabf} shows that $(\partial_j f_n)$ is Cauchy in $L^2(\B^5_R)$ for any $j\in \{1,2,\dots,5\}$.
We set $g_j:=\lim_{n\to\infty}\partial_j f_n \in L^2(\B^5_R)$.
By dominated convergence we infer
\[ \int_{\B^5_R}g_j \varphi=\lim_{n\to\infty}\int_{\B^5_R}\partial_j f_n \varphi 
=-\lim_{n\to\infty}\int_{\B^5_R}f_n\partial_j\varphi=-\int_{\B^5_R}f\partial_j \varphi \]
for any $\varphi\in C^\infty_c(\B^5_R)$.
Consequently, $\partial_j f=g_j$ in the weak sense and this shows $\nabla f\in L^2(\B^5_R)$ with the bound \eqref{eq:nabf}.

Let $f\in C^\infty_c(\R^5)$. Then we have
\begin{align*}
\|\mc L f\|_{L^2(\B_R^5)}^2&=(\Delta f-\Lambda f+V_0f | \Delta f-\Lambda f+V_0 f)_{L^2(\B_R^5)} \\
&=\|\Delta f\|_{L^2(\B_R^5)}^2+2(\Delta f | -\Lambda f+V_0 f)_{L^2(\B_R^5)}
+\|-\Lambda f+V_0 f\|_{L^2(\B_R^5)}^2
\end{align*}
which yields the bound
\begin{align*} 
\|\Delta f\|_{L^2(\B_R^5)}&\lesssim \|\mc L f\|_{L^2(\B_R^5)} 
+\|\Lambda f\|_{L^2(\B^5_R)}+\|V_0 f\|_{L^2(\B_R^5)} \\
&\leq C_R\|f\|_{G(\mc L)} 
+C_R \|\nabla f\|_{L^2(\B_R^5)} \\
&\leq C_R\|f\|_{G(\mc L)} 
\end{align*}
by Eq.~\eqref{eq:nabf}.
Consequently, a density argument as above finishes the proof.
\end{proof}

In order to control the full Sobolev norm
\[ \|f\|_{H^k(\B^5_R)}=\sum_{|\alpha|\leq k}\|\partial^\alpha f\|_{L^2(\B_R^5)} \]
for $k=2$,
we need two technical results which are completely elementary since we restrict ourselves to radial functions.
First, we have a trace lemma.

\begin{lemma}
\label{lem:trace}
Let $R\geq 1$. Then we have the bound
\[ \|\nabla f\|_{L^\infty(\partial\B^5_R)}\leq C_R \|\Delta f\|_{L^2(\B^5_R)} \]
for all radial $f\in C^2(\overline{\B^5_R})$.
\end{lemma}

\begin{proof}
By assumption, there exists a function $\tilde f\in C^2([0,R])$ such that $f(x)=\tilde f(|x|)$.
The fundamental theorem of calculus yields
\begin{align*}
|R^4 \tilde f'(R)|=\left |\int_0^R \partial_r [r^4 \tilde f'(r)]\d r\right |
\leq C_R \left (\int_0^R \left |\tfrac{1}{r^4}\partial_r [r^4 \tilde f'(r)]\right |^2 r^8 \d r \right )^{1/2}\leq C_R\|\Delta f\|_{L^2(\B^5_R)}
\end{align*}
and thus,
\[ |\partial_j f(x)|=\left |\tfrac{x_j}{|x|}\tilde f'(|x|)\right |\leq C_R \|\Delta f\|_{L^2(\B^5_R)} \]
for all $x\in \partial\B^5_R$ and $j\in \{1,2,\dots,5\}$.
\end{proof}

Next, by an extension argument and Fourier analysis, we easily get control on mixed derivatives.
Here and in the following, $\mc F$ is the Fourier transform
\[ (\mc F f)(\xi):=\int_{\R^5}e^{-\I \xi x}f(x)\d x. \]

\begin{lemma}
\label{lem:mixed}
Let $R\geq 1$. Then we have the bound
\[ \|\partial_j\partial_k f\|_{L^2(\B^5_R)}\leq C_R\left (\|\Delta f\|_{L^2(\B^5_R)}+ \|\nabla f\|_{L^2(\B^5_R)}+\|f\|_{L^2(\B^5_R)}\right ) \]
for all radial $f\in C^2(\overline{\B^5_R})$  and all $j,k\in \{1,2,\dots,5\}$.
\end{lemma}

\begin{proof}
Let $\tilde f\in C^2([0,R])$ such that $f(x)=\tilde f(|x|)$. We define an extension $\tilde{\mc E}\tilde f$ of $\tilde f$ by
\[ \tilde{\mc E} \tilde f(r):=\left \{\begin{array}{ll}\tilde f(r) & \mbox{ for }r\in [0,R] \\
\tilde f(2R-r)-2\tilde f'(R)(R-r) & \mbox{ for }r\in (R,2R)
\end{array} \right . .
\]
Then we have
\begin{align*}
 \lim_{r\to R+}\tilde{\mc E}\tilde f(r)&=\tilde f(R) \\
 \lim_{r\to R+}(\tilde{\mc E}\tilde f)'(r)&=\lim_{r\to R+}\big [-\tilde f'(2R-r)+2\tilde f'(R)\big ]=\tilde f'(R) \\
 \lim_{r\to R+}(\tilde{\mc E}\tilde f)''(r)&=\lim_{r\to R+}\tilde f''(2R-r)=\tilde f''(R)
 \end{align*}
 and thus, $\tilde{\mc E}\tilde f\in C^2([0,2R))$.
 Furthermore,
 \begin{align*}
  \int_R^{\frac32 R} |\tilde{\mc E}\tilde f(r)|^2 r^4 \d r
  &=\int_{\frac12 R}^R |\tilde {\mc E}\tilde f(2R-r)|^2(2R-r)^4 \d r
  \leq C_R \int_{\frac12 R}^R |\tilde {\mc E}\tilde f(2R-r)|^2 r^4 \d r \\
  &\leq C_R \int_{\frac12 R}^R |\tilde f(r)|^2 r^4 \d r+C_R |\tilde f'(R)|^2.
  \end{align*}
  Analogously, we obtain
  \begin{equation}
  \label{eq:ext}
  \||\cdot|^2(\tilde {\mc E} \tilde f)^{(k)}\|_{L^2(R,\frac32 R)}\leq C_R 
  \||\cdot|^2 \tilde f^{(k)}\|_{L^2(\frac12R,R)}+C_R |\tilde f'(R)|
  \end{equation}
  for any $k\in \{0,1,2\}$.
  
Now let $\chi: \R^5 \to [0,1]$ be a smooth cut-off that satisfies $\chi(x)=1$ for $|x|\leq 1$ and $\chi(x)=0$ for $|x|\geq \frac32$ and set 
\[ \mc E f(x):=\chi(\tfrac{x}{R})\tilde{\mc E}\tilde f(|x|). \]
Then $\mc E f\in C_c^2(\R^5)$ with $\supp(\mc E f)\subset \overline{\B^5_{\frac32R}}$ and 
$\mc Ef=f$ on $\B^5_R$.
From Eq.~\eqref{eq:ext} and Lemma \ref{lem:trace} we obtain the bound
\begin{align*} \|\Delta \mc E f\|_{L^2(\R^5)}&\simeq \|\Delta \mc E f\|_{L^2(\B^5_R)}
+\|\Delta \mc E f\|_{L^2(\B^5_{\frac32R}\setminus \B^5_R)} \\
&\leq C_R \left (\|\Delta f\|_{L^2(\B^5_R)}+\|\nabla f\|_{L^2(\B^5_R)}+\|f\|_{L^2(\B^5_R)}\right ) .
\end{align*}
Consequently, the estimate
\begin{align*} 
\|\partial_j\partial_k f\|_{L^2(\B^5_R)}\leq \|\partial_j\partial_k \mc E f\|_{L^2(\R^5)} 
\lesssim \||\cdot|^2 \mc F\mc E f\|_{L^2(\R^5)}
\simeq \|\Delta \mc E f\|_{L^2(\R^5)}
\end{align*}
finishes the proof.
\end{proof}

Now we can control the linear evolution on the Sobolev space $H^2(\B^5_R)$.

\begin{corollary}
\label{cor:H2B}
Let $f\in \mc D(\mc L)$ and $R\geq 1$. Then $e^{s\mc L}f \in H^2(\B^5_R)$ for all $s\geq 0$ and we have the bound
\[ \|e^{s\mc L}f\|_{H^2(\B^5_R)}
\leq C_R e^{-c_0 s}\|f\|_{G(\mc L)} \]
for all $s\geq 0$ and all $f\in \mc D(\mc L)$ satisfying $(f|\psi_1)_{L^2_\sigma(\R^5)}=0$.
\end{corollary}

\begin{proof}
By Lemma \ref{lem:mixed} it suffices to control $\nabla e^{s\mc L}f$ and $\Delta e^{s\mc L}f$.
Since $\mc D(\mc L)$ is invariant under $e^{s\mc L}$, Lemma \ref{lem:H2B} implies $\nabla e^{s\mc L}f, \Delta e^{s\mc L}f \in L^2(\B^5_R)$. Consequently, from Lemma \ref{lem:H2B}
and Proposition \ref{prop:gen} we infer
\begin{align*}
 \|\Delta e^{s\mc L}f\|_{L^2(\B^5_R)}+\|\nabla e^{s\mc L}f\|_{L^2(\B^5_R)}
 &\leq C_R \left (\|\mc L e^{s\mc L}f\|_{L^2_\sigma(\R^5)}+\|e^{s\mc L}f\|_{L^2_\sigma(\R^5)} \right ) \\
 &=C_R \left (\|e^{s\mc L}\mc L f\|_{L^2_\sigma(\R^5)}+\|e^{s\mc L}f\|_{L^2_\sigma(\R^5)} \right ) \\
 &\leq C_R e^{-c_0 s}\|f\|_{G(\mc L)}
 \end{align*}
 since
 \[ (\mc L f|\psi_1)_{L^2_\sigma(\R^5)}=(f|\mc L\psi_1)_{L^2_\sigma(\R^5)}
 =(f|\psi_1)_{L^2_\sigma(\R^5)}=0. \]
\end{proof}

Next, we improve the above by two derivatives.

\begin{lemma}
\label{lem:H4B}
Let $f\in \mc D(\mc L^2)$ and $R\geq 1$.
Then $\nabla \Delta f, \Delta^2 f\in L^2(\B^5_R)$ and 
we have the bound 
\[ \|\Delta^2 f\|_{L^2(\B_R^5)}+\|\nabla \Delta f\|_{L^2(\B^5_R)}\leq C_R \|f\|_{G(\mc L^2)} \]
for all $R\geq 1$ and all $f\in \mc D(\mc L^2)$.
\end{lemma}

\begin{proof}
Let $f\in C^\infty_c(\R^5)$ and $R\geq 1$.
From Lemma \ref{lem:H2B} we have the bound
\[ \|\nabla \mc L f\|_{L^2(\B^5_R)}\leq C_R \|f\|_{G(\mc L^2)}. \]
Expanding the square yields
\begin{align*} \|\nabla \mc L f\|_{L^2(\B^5_R)}^2&=\|\nabla(\Delta-\Lambda+V_0)f\|_{L^2(\B^5_R)}^2 \\
&=\|\nabla\Delta f\|_{L^2(\B^5_R)}^2+2(\nabla\Delta f | \nabla(-\Lambda+V_0)f)_{L^2(\B^5_R)} 
+\|\nabla(-\Lambda+V_0)f\|_{L^2(\B^5_R)}^2
\end{align*}
and thus,
\begin{align*} \|\nabla\Delta f\|_{L^2(\B^5_R)}&\lesssim \|\nabla\mc Lf\|_{L^2(\B^5_R)}
+\|f\|_{H^2(\B^5_R)}\leq C_R \|f\|_{G(\mc L^2)}
\end{align*}
by Lemmas \ref{lem:H2B} and \ref{lem:mixed}.

For $\Delta^2 f$ we expand $\|\Delta \mc L f\|_{L^2(\B^5_R)}^2$ and use Lemma \ref{lem:H2B} together with the bound on $\nabla\Delta f$ to obtain
\[ \|\Delta^2 f\|_{L^2(\B^5_R)} \lesssim \|\Delta \mc L f\|_{L^2(\B^5_R)}^2
+\|\nabla\Delta f\|_{L^2(\B^5_R)}^2
+C_R\|f\|_{G(\mc L)}^2\leq C_R \|f\|_{G(\mc L^2)}^2. \]
\end{proof}

\begin{corollary}
\label{cor:H4B}
Let $f\in \mc D(\mc L^2)$ and $R\geq 1$. Then $e^{s\mc L}f \in H^4(\B^5_R)$ for all $s\geq 0$ and we have the bound
\[ \|e^{s\mc L}f\|_{H^4(\B^5_R)}
\leq C_R e^{-c_0 s}\|f\|_{G(\mc L^2)} \]
for all $s\geq 0$ and all $f\in \mc D(\mc L^2)$ satisfying $(f|\psi_1)_{L^2_\sigma(\R^5)}=0$.
\end{corollary}

\begin{proof}
It suffices to note that $\mc D(\mc L^2)$ is invariant under $e^{s\mc L}$ so that Lemma \ref{lem:H4B} can be applied to $e^{s\mc L}f$. Proposition \ref{prop:gen} and Lemma \ref{lem:mixed} then yields the  statement.
\end{proof}

\subsection{Estimates in unweighted global Sobolev norms}

Next, we prove bounds in $\dot H^2(\R^5)$ and $\dot H^4(\R^5)$.
The intersection $\dot H^2(\R^5)\cap \dot H^4(\R^5)$ will be our main space where we study the evolution.
First, we have to ensure that unweighted Sobolev spaces are invariant under $e^{s\mc L}$.

\begin{lemma}
\label{lem:Hk}
Let $k\in \N_0$ and $f\in H^k(\R^5)$.
Then we have $e^{s\mc L}f\in H^k(\R^5)$ for all $s\geq 0$ and $e^{s\mc L}$ is a strongly continuous semigroup on $H^k(\R^5)$.
\end{lemma}

\begin{proof}
We denote by $\mc L_0=\mc L-V_0+\frac12$ the principal part of $\mc L$.
The operator $\mc L_0: \mc D(\mc L)\subset L^2_\sigma(\R^5)\to L^2_\sigma(\R^5)$ is self-adjoint and it generates the semigroup $e^{s\mc L_0}$.
As a matter of fact, $e^{s\mc L_0}$ can be given explicitly and we have

\[ [e^{s\mc L_0}f](x)=(K_s * f)(e^{-s/2} x) \]
where
\[ K_s(x)=[\pi \alpha(s)]^{-5/2}e^{-|x|^2/\alpha(s)},\qquad \alpha(s)=4(1-e^{-s}). \]
This is easily verified by an explicit computation.
Since $K_s\in L^1(\R^5)$ for any $s> 0$, dominated convergence and Young's inequality immediately imply the invariance of $H^k(\R^5)$ under $e^{s\mc L_0}$.
By rescaling we infer
\begin{align*}
 e^{s\mc L_0}f(x)-f(x)&=[\pi\alpha(s)]^{-\frac52}\int_{\R^5}
e^{-|x'|^2/\alpha(s)}\big [f(e^{-s/2}x-x')-f(x)\big ]\d x' \\
&=\pi^{-\frac52}\int_{\R^5}\big [f(e^{-s/2}x-\alpha(s)^\frac12 x')-f(x)\big ]
e^{-|x'|^2}\d x'
\end{align*}
and Minkowski's inequality yields
\[ \|e^{s\mc L_0}f-f\|_{L^2(\R^5)}\lesssim \int_{\R^5}\big \|f(e^{-s/2}(\cdot)-\alpha(s)^\frac12 x')-f\big \|_{L^2(\R^5)}e^{-|x'|^2}\d x'. \]
Since scaling and translation are continuous operations on $L^2(\R^5)$,
we infer
\[ \|f(e^{-s/2}(\cdot)-\alpha(s)^\frac12 x')-f\big \|_{L^2(\R^5)}\to 0 \]
as $s\to 0+$ for any fixed $x'\in \R^5$.
Consequently, by dominated convergence, we obtain
\[ \|e^{s\mc L_0}f-f\|_{L^2(\R^5)}\to 0 \]
as $s\to 0+$.
The same argument yields
$\|e^{s\mc L_0}f-f\|_{H^k(\R^5)}\to 0$
as $s\to 0+$.
We conclude that $e^{s\mc L_0}$ is strongly continuous on $H^k(\R^5)$.
Evidently, the map $f\mapsto V_0 f$ is bounded on $H^k(\R^5)$ and thus, by the bounded perturbation theorem, $e^{s\mc L}$ is a strongly continuous semigroup on $H^k(\R^5)$.
\end{proof}

\begin{lemma}
\label{lem:dotH2}
Let $f\in  \mc D(\mc L)\cap \dot H^2(\R^5)$. Then $e^{s\mc L}f\in \dot H^2(\R^5)$ for all $s\geq 0$ and there exists a constant $c_1>0$ such that
\[ \|e^{s\mc L}f\|_{\dot H^2(\R^5)}\lesssim e^{-c_1 s}\left (\|f\|_{G(\mc L)}+\|f\|_{\dot H^2(\R^5)} \right ) \]
for all $s\geq 0$ and all $f\in \mc D(\mc L)\cap \dot H^2(\R^5)$ satisfying $(f|\psi_1)_{L^2_\sigma(\R^5)}=0$.
\end{lemma}

\begin{proof}
Let $f\in C^\infty_c(\R^5)$ and
note the commutator relation $[\Delta, \Lambda]f=\Delta f$ which yields
\begin{align*} 
\Delta \mc L f&=\Delta^2 f-\Delta \Lambda f+\Delta(V_0 f)
=\Delta^2 f-\Lambda \Delta f-\Delta f+\Delta(V_0 f).
\end{align*}
Consequently, with $(-\Lambda \Delta f|\Delta f)_{L^2(\R^5)}= \frac34 \|\Delta f\|_{L^2(\R^5)}^2$ we obtain
\begin{align}
\label{eq:DLf}
\begin{split}
(\Delta \mc L f | \Delta f)_{L^2(\R^5)}&=-\|\nabla\Delta f\|_{L^2(\R^5)}^2
-(\Lambda \Delta f |\Delta f)_{L^2(\R^5)}-\|\Delta f\|_{L^2(\R^5)}^2
+(\Delta(V_0 f)|\Delta f)_{L^2(\R^5)} \\
&\leq -\|\nabla\Delta f\|_{L^2(\R^5)}^2-\tfrac14 \|\Delta f\|_{L^2(\R^5)}^2
+(\Delta(V_0 f)|\Delta f)_{L^2(\R^5)}.
\end{split}
\end{align}
Now we claim the estimate 
\begin{equation}
\label{eq:mixed}
|(\Delta(V_0 f)|\Delta f)_{L^2(\R^5)}|\leq C_R \|f\|_{G(\mc L)}^2+
\tfrac{C}{R^2}\|\Delta f\|_{L^2(\R^5)}^2
\end{equation}
for all $R\geq 1$.
To prove this, we note that $\Delta(V_0 f)=\Delta V_0 f+2\nabla V_0 \nabla f+V_0\Delta f$ and estimate each of these terms individually.
Clearly, 
\begin{align*} |(V_0 \Delta f | \Delta f)_{L^2(\R^5)}|&\lesssim \||V_0|^\frac12 \Delta f\|_{L^2(\R^5)}^2
= \||V_0|^\frac12 \Delta f\|_{L^2(\B^5_R)}^2+
\||V_0|^\frac12 \Delta f\|_{L^2(\R^5\setminus \B^5_R)}^2 \\
&\leq C_R \|f\|_{G(\mc L)}^2+\tfrac{C}{R^2}\|\Delta f\|_{L^2(\R^5)}^2
\end{align*}
where we have used Lemma \ref{lem:H2B} and the decay $|V_0(x)|\lesssim \langle x\rangle^{-2}$.
Next,
\begin{align*}
|(\nabla V_0 \nabla f | \Delta f)_{L^2(\R^5)}|&\lesssim
\||\nabla V_0|^\frac23 \nabla f\|_{L^2(\R^5)}^2+\||\nabla V_0|^\frac13 \Delta f\|_{L^2(\R^5)}^2.
\end{align*}
Thanks to the decay $|\nabla V_0(x)|\lesssim \langle x\rangle^{-3}$, the last term can be estimated as before.
For the first term we use the decay of $\nabla V_0$, Lemma \ref{lem:H2B}, and Hardy's inequality to estimate
\begin{align*}
\||\nabla V_0|^\frac23 \nabla f\|_{L^2(\R^5)}&\lesssim
\|\langle\cdot\rangle^{-2}\nabla f\|_{L^2(\R^5)}
\simeq \|\nabla f\|_{L^2(\B^5_R)}+\||\cdot|^{-2}\nabla f\|_{L^2(\R^5\setminus \B^5_R)} \\
&\leq C_R \|f\|_{G(\mc L)}+\tfrac{C}{R}\||\cdot|^{-1}\nabla f\|_{L^2(\R^5\setminus \B^5_R)} \\
&\leq C_R \|f\|_{G(\mc L)}+\tfrac{C}{R}\||\cdot|^{-1}\nabla f\|_{L^2(\R^5)}  \\
&\leq C_R \|f\|_{G(\mc L)}+\tfrac{C}{R}\|\Delta f\|_{L^2(\R^5)}.
\end{align*}
In view of the decay $|\Delta V_0(x)|\lesssim \langle x \rangle^{-4}$, the term $(\Delta V_0 f|\Delta f)_{L^2(\R^5)}$ can be estimated analogously.
This proves Eq.~\eqref{eq:mixed}.

Having Eq.~\eqref{eq:mixed} at our disposal, we obtain from Eq.~\eqref{eq:DLf} the bound
\begin{equation}
\label{eq:dissdotH2}
 (\Delta \mc L f | \Delta f)_{L^2(\R^5)}\leq (-\tfrac14+\tfrac{C}{R^2})\|\Delta f\|_{L^2(\R^5)}^2+C_R \|f\|_{G(\mc L)}^2. 
\end{equation}
By approximation, Eq.~\eqref{eq:dissdotH2} extends to all $f\in \mc D (\mc L)$ satisfying $\mc L f, f \in H^2(\R^5)$.
From Lemma \ref{lem:Hk} we know that $\mc L e^{s\mc L}f, e^{s\mc L}f\in H^2(\R^5)$ and Eq.~\eqref{eq:dissdotH2} yields
\begin{align*}
\tfrac12 \partial_s \|\Delta e^{s\mc L}f\|_{L^2(\R^5)}^2&=
(\partial_s \Delta e^{s\mc L}f | \Delta e^{s\mc L}f )_{L^2(\R^5)}=(\Delta \mc L e^{s\mc L}f |
\Delta e^{s\mc L}f)_{L^2(\R^5)} \\
&\leq (-\tfrac14+\tfrac{C}{R^2})\|\Delta e^{s\mc L}f\|_{L^2(\R^5)}^2
+C_R \|e^{s\mc L}f\|_{G(\mc L)}^2 \\
&\leq -\tfrac18 \|\Delta e^{s\mc L}f\|_{L^2(\R^5)}^2 + C_R e^{-2c_0 s}\|f\|_{G(\mc L)}^2
\end{align*}
by choosing $R\geq 1$ sufficiently large.
From now on $R$ is fixed and hence, $C_R=C$.
Upon setting $c_1=\frac12\min\{c_0,\frac18\}>0$, we infer
\[ \tfrac12 \partial_s \|\Delta e^{s\mc L}f\|_{L^2(\R^5)}^2\leq -2c_1 
\|\Delta e^{s\mc L}f\|_{L^2(\R^5)}^2 + C e^{-4c_1 s}\|f\|_{G(\mc L)}^2 \]
and this inequality may be rewritten as
\[ \tfrac12 \partial_s \left [e^{4c_1 s}\|\Delta e^{s\mc L}f\|_{L^2(\R^5)}^2 \right ]
\leq C \|f\|_{G(\mc L)}^2. \]
Consequently, integration yields the bound
\begin{align*}
 \|\Delta e^{s\mc L}f\|_{L^2(\R^5)}^2&\lesssim \langle s\rangle e^{-4c_1 s} \left (\|f\|_{G(\mc L)}^2+\|f\|_{\dot H^2(\R^5)}^2 \right )  \\
 &\lesssim e^{-2c_1 s} \left (\|f\|_{G(\mc L)}^2+\|f\|_{\dot H^2(\R^5)}^2 \right ) .
 \end{align*}
 By a density argument, this bound holds for all $f\in \mc D(\mc L)\cap \dot H^2(\R^5)$.
\end{proof}

It is now straightforward to upgrade to $\dot H^4$.

\begin{lemma}
\label{lem:dotH4}
Let $f\in  \mc D(\mc L^2)\cap \dot H^2(\R^5)\cap \dot H^4(\R^5)$. Then $e^{s\mc L}f\in \dot H^4(\R^5)$ for all $s\geq 0$ and there exists a constant $c_1>0$ such that
\[ \|e^{s\mc L}f\|_{\dot H^4(\R^5)}\lesssim e^{-c_1 s}\left (\|f\|_{G(\mc L^2)}+\|f\|_{\dot H^2(\R^5)}+\|f\|_{\dot H^4(\R^5)} \right ) \]
for all $s\geq 0$ and all $f\in \mc D(\mc L^2)\cap \dot H^2(\R^5)\cap \dot H^4(\R^5)$ satisfying $(f|\psi_1)_{L^2_\sigma(\R^5)}=0$.
\end{lemma}

\begin{proof}
Let $f\in C^\infty_c(\R^5)$. 
By applying the commutator relation $[\Delta,\Lambda]f=\Delta f$ twice, we obtain the estimate
\[ (\Delta^2\mc Lf |\Delta^2 f)_{L^2(\R^5)}\leq -\tfrac54 \|\Delta^2 f\|_{L^2(\R^5)}^2
+(\Delta^2 (V_0 f)|\Delta^2 f)_{L^2(\R^5)}, \]
cf.~Eq.~\eqref{eq:DLf}.
Consequently, it suffices to follow the logic in the proof of Lemma \ref{lem:dotH2} and apply Lemma \ref{lem:H4B}.
\end{proof}

\subsection{Control of the linearized flow}

Finally, we arrive at the main result on the linearized flow. First, we define the main Sobolev space we will be working with and prove an elementary embedding result.

\begin{definition}
The Banach space $X$ is defined as the completion of all radial functions in $C^\infty_c(\R^5)$ with respect to the norm
\[ \|f\|_X=\|f\|_{\dot H^2(\R^5)\cap \dot H^4(\R^5)}=\|\Delta f\|_{L^2(\R^5)}+\|\Delta^2 f\|_{L^2(\R^5)}. \] 
\end{definition}

\begin{lemma}
\label{lem:Linf}
Let $s\in [0,\frac32)$. Then we have the bound
\[ \||\nabla|^s f\|_{L^\infty(\R^5)}\lesssim \|f\|_X \]
for all $f\in C^\infty_c(\R^5)$.
\end{lemma}

\begin{proof}
We readily estimate
\begin{align*}
\||\nabla|^s f\|_{L^\infty(\R^5)}&\lesssim \||\cdot|^s \mc F f\|_{L^1(\R^5)}\simeq
\||\cdot|^s \mc F f\|_{L^1(\B^5)}+\||\cdot|^s\mc F f\|_{L^1(\R^5\setminus \B^5)} \\
&\lesssim \||\cdot|^{-2}\|_{L^2(\B^5)}\||\cdot|^2 \mc F f\|_{L^2(\R^5)}
+\||\cdot|^{-4+s}\|_{L^2(\R^5\setminus \B^5)}\||\cdot|^4\mc F f\|_{L^2(\R^5)} \\
&\lesssim \|f\|_{\dot H^2(\R^5)}+\|f\|_{\dot H^4(\R^5)}.
\end{align*}
\end{proof}

Now we can prove the following simple but useful embedding theorem.

\begin{lemma}
\label{lem:embed}
We have the continuous embeddings
\[ H^4_{\mathrm{rad}}(\R^5)\hookrightarrow X \hookrightarrow C^1(\R^5)\cap W^{1,\infty}(\R^5) \]
where $H^4_\mathrm{rad}(\R^5)=\{f\in H^4(\R^5): f \mbox{ radial}\}$.
\end{lemma}

\begin{proof}
Let $f\in H^4_{\mathrm{rad}}(\R^5)$. Then there exists a sequence $(f_n)_{n\in \N}\subset C^\infty_c(\R^5)$ of radial functions such that $f_n\to f$ with respect to $\|\cdot\|_{H^4(\R^5)}$.
This implies that $(f_n)_{n\in\N}$ is Cauchy with respect to $\|\cdot\|_X$ and thus, there exists a limiting element $\hat f \in X$ such that $f_n\to \hat f$ in $X$.
We define a map $\iota: H^4_\mathrm{rad}(\R^5)\to X$ by setting $\iota(f):=\hat f$.
Obviously, $\iota$ is linear.
We claim that $\iota$ is injective. Indeed, if $\iota(f)=0$, there exists a sequence $(f_n)_{n\in \N}\subset C^\infty_c(\R^5)$ that converges to $f$ in $H^4_{\mathrm{rad}}(\R^5)$ and to $0$ in $X$.
By Lemma \ref{lem:Linf} we see that $\lim_{n\to\infty}\|f_n\|_{L^\infty(\R^5)}=0$. In particular, $f_n\wto 0$ in $L^2(\R^5)$.
On the other hand, $f_n\to f$ in $H^4_\mathrm{rad}(\R^5)$ implies $f_n\wto f$ in $L^2(\R^5)$ and the uniqueness of weak limits shows that $f=0$.
Clearly, we have $\|\iota(f)\|_X\lesssim \|f\|_{H^4(\R^5)}$ and thus, $\iota: H^4_\mathrm{rad}(\R^5)\to X$ is a continuous embedding.

The second assertion is proved similarly. Indeed, given $f\in X$ we find a sequence $(f_n)_{n\in\N}\subset C^\infty_c(\R^5)$ such that $f_n\to f$ in $X$.
By Lemma \ref{lem:Linf}, $(f_n)_{n\in \N}$ is Cauchy in $W^{1,\infty}(\R^5)$ and therefore converges to a limiting function $\hat f\in C^1(\R^5)\cap W^{1,\infty}(\R^5)$.
Using this, we define an inclusion map $\iota: X\to C^1(\R^5)\cap W^{1,\infty}(\R^5)$ by setting $\iota(f):=\hat f$.
It remains to show that $\iota$ is injective. If $\iota(f)=0$, it follows that there exists a sequence $(f_n)_{n\in\N}\subset C^\infty_c(\R^5)$ that converges to $f$ in $X$ and to $0$ in $L^\infty(\R^5)$. Consequently,
\[ \left |\int_{\R^5} \Delta f_n \varphi\right | = \left | \int_{\R^5} f_n \Delta \varphi
\right |\lesssim \|f_n\|_{L^\infty(\R^5)}\to 0 \]
for any $\varphi\in C_c^\infty(\R^5)$ and thus,
$\Delta f_n \wto 0$ in $L^2(\R^5)$.
Analogously, we obtain $\Delta^2 f_n\wto 0$ in $L^2(\R^5)$.
By the uniqueness of weak limits we therefore have
$\lim_{n\to\infty}\|f_n\|_X=0$
and this shows $f=0$.
\end{proof}

\begin{theorem}
\label{thm:lin}
The Sobolev space $X$ is invariant under $e^{s\mc L}$ and there exists a constant $\omega_0>0$ such that
\[ \|e^{s\mc L}f\|_X\lesssim e^{-\omega_0 s}\|f\|_X \]
for all $s\geq 0$ and all $f\in X$ satisfying $(f|\psi_1)_{L^2_\sigma(\R^5)}=0$.
\end{theorem}

\begin{proof}
By Lemma \ref{lem:embed}, $\mc D(\mc L^2)\cap H^4(\R^5)\hookrightarrow X$. Since the former space is invariant under $e^{s\mc L}$, see Lemma \ref{lem:Hk}, it follows that $e^{s\mc L} f\in X$ for all $s\geq 0$ and all $f\in C^\infty_c(\R^5)$.
Consequently, in view of Lemmas \ref{lem:dotH2}, \ref{lem:dotH4}, and a density argument, it suffices to prove the bound
\[ \|f\|_{G(\mc L^2)}\lesssim \|f\|_X \]
for all $f\in C^\infty_c(\R^5)$.
Thanks to the strong decay of the weight $\sigma(x)=e^{-|x|^2/4}$, we immediately obtain
\begin{align*}
\|f\|_{G(\mc L^2)}&\lesssim \||\cdot|^{-2}f\|_{L^2(\R^5)}+\||\cdot|^{-1} \nabla f\|_{L^2(\R^5)}+\|\Delta f\|_{L^2(\R^5)} \\
&\quad +\||\cdot|^{-1}\nabla \Delta f\|_{L^2(\R^5)}
+\|\Delta^2 f\|_{L^2(\R^5)} \\
&\lesssim \|f\|_X
\end{align*}
by Hardy's inequality.
\end{proof}

\section{The nonlinear evolution}

\noindent Now we turn to the full nonlinear problem Eq.~\eqref{eq:full}.
As before with the linear operator, we switch to $5$-dimensional notation and define the nonlinearity $\mc N$, acting on functions $f:\R^5\to\R$, by
\[ \mc N(f)(x):=-\frac{1}{|x|^3}\left [\sin(2f_0(|x|)+2|x|f(x))-\sin(2f_0(|x|))-2|x|\cos(2f_0(|x|))f(x)\right ]. \]
With this convention, Eq.~\eqref{eq:full} can be written as
\begin{equation}
\label{eq:fullop}
\left \{ \begin{array}{l}
\partial_s \phi(s)=\mc L\phi(s)+\mc N(\phi(s)) \\
\phi(0)=\mc U(h,T)
\end{array} \right .
\end{equation}
where 
\[ \mc U(h,T)(x):=f_0(\sqrt T |x|)/|x|-f_0(|x|)/|x|+h(\sqrt T |x|)/|x| \]
and $\phi: [0,\infty)\to X$.

So far, this is purely formal. In what follows we first prove basic embedding theorems and then some Moser-type inequalities. These will allow us to show that the nonlinearity is locally Lipschitz on $X$. Next, we study mapping properties of the ``initial data operator'' $\mc U$ and finally, we implement an infinite-dimensional version of the Lyapunov-Perron method to prove global existence for Eq.~\eqref{eq:fullop}.

\subsection{Further properties of the space $X$}

\begin{corollary}[Algebra property]
\label{cor:alg}
We have the bound
\[ \|fg\|_X\lesssim \|f\|_X \|g\|_X \]
for all $f,g\in X$.
As a consequence, $X$ is a Banach algebra.
\end{corollary}

\begin{proof}
This is a straightforward consequence of the Leibniz rule, the Gagliardo-Nirenberg inequality (see e.g.~\cite{Tay11}), and Lemma \ref{lem:Linf}.
\end{proof}

Next, we prove weighted $L^\infty$ bounds outside of balls.
As opposed to Lemma \ref{lem:Linf} and Corollary \ref{cor:alg}, the restriction to radial functions is crucial here.

\begin{lemma}
\label{lem:LinfB}
We have the bounds
\begin{align*}
\||\cdot|^\frac32 f\|_{L^\infty(\R^5\setminus \B^5)}&\lesssim \|f\|_{\dot H^1(\R^5)} \\
 \||\cdot|^\frac12 f\|_{L^\infty(\R^5\setminus \B^5)}&\lesssim \|f\|_{\dot H^2(\R^5)}  
 \end{align*}
 for all \emph{radial} $f\in C^\infty_c(\R^5)$.
\end{lemma}

\begin{proof}
Let $f\in C^\infty_c(\R^5)$ be radial and write $f(x)=\tilde f(|x|)$.
The fundamental theorem of calculus yields
\[ \tilde f(r)=-\int_r^\infty \tilde f'(s)\d s=-\int_r^\infty s^{-2}\tilde f'(s)s^2 \d s \]
and thus, by Cauchy-Schwarz,
\[ |\tilde f(r)|\leq \||\cdot|^2 \tilde f'\|_{L^2(1,\infty)}\left (
\int_r^\infty s^{-4}\d s\right )^{1/2}\lesssim r^{-\frac32} \||\cdot|^2 \tilde f'\|_{L^2(1,\infty)}\lesssim r^{-\frac32}\|\nabla f\|_{L^2(\R^5)} \]
for all $r\geq 1$. This implies the first assertion.

For the second statement we proceed similarly and use
\[ \tilde f(r)=\int_r^\infty \int_s^\infty \tilde f''(t)\d t =
\int_r^\infty \int_s^\infty t^{-2}\tilde f''(t)t^2\d t \]
to obtain the bound
\[ |\tilde f(r)|\leq \||\cdot|^2 \tilde f''\|_{L^2(1,\infty)}
\int_r^\infty \left (\int_s^\infty t^{-4}\d t \right )^{1/2}\d s
\lesssim r^{-\frac12}\||\cdot|^2 \tilde f''\|_{L^2(1,\infty)} \]
for all $r\geq 1$.
Now note that
\[ \tilde f'(|x|)=\tfrac{x^j}{|x|}\partial_j f(x) \]
and thus, 
by Hardy's inequality, we infer
\[ r^\frac12 |\tilde f(r)|\lesssim \|\Delta f\|_{L^2(\R^5)}+\||\cdot|^{-1}\nabla f\|_{L^2(\R^5)} 
\lesssim \|\Delta f\|_{L^2(\R^5)} \]
for all $r\geq 1$, which is the desired result.

\end{proof}

\subsection{Nonlinear estimates}

For $\delta>0$ we set
\[ X_\delta:=\{f\in X: \|f\|_X\leq \delta\}. \]
The goal of this section is to prove that the nonlinearity $\mc N$ is locally Lipschitz on $X$.
The key results in this respect are the following Moser-type inequalities.
First, we focus on large radii where we need to assume a decay property.

\begin{proposition}
\label{prop:MoseroB}
Let $\Phi\in C^4(\R\times \R^5)$ and suppose 
\[ |\partial^\alpha \Phi(v,x)|\lesssim \langle x\rangle^{-1} \]
for all $(v,x)\in \R\times \R^5\setminus\B^5$ and all multi-indices $\alpha\in \N_0^6$ with $|\alpha|\leq 4$. For $f: \R^5\to\R$ set
\[ \mc M(f)(x):=f(x)^2 \Phi(|x|f(x),x). \]
Then we have the bound
\[ \|\mc M(f)-\mc M(g)\|_{\dot H^2(\R^5\setminus \B^5)\cap \dot H^4(\R^5\setminus \B^5)}\lesssim (\|f\|_X+\|g\|_X)\|f-g\|_X \]
for all $f,g \in X_1 \cap C^\infty_c(\R^5)$.
\end{proposition}

\begin{proof}
Let $f,g \in X_1\cap C^\infty_c(\R^5)$ and set $\mc I(f)(x):=\Phi(|x|f(x),x)$. Then we have
\[ \mc M(f)-\mc M(g)=f^2 \mc I(f)-g^2\mc I(g)=(f^2-g^2)\mc I(f)+g^2[\mc I(f)-\mc I(g)] \]
and
\begin{align*}
\mc I(f)(x)-\mc I(g)(x)&=\int_0^1 \partial_t \Phi\big (|x|g(x)+t|x|(f(x)-g(x)),x\big )\d t \\
&=[f(x)-g(x)]|x|\int_0^1 \partial_1 \Phi\big (|x|g(x)+t|x|(f(x)-g(x)),x\big )\d t \\
&=:[f(x)-g(x)]\mc J(f,g)(x).
\end{align*}
Consequently, it suffices to prove
\begin{align*}
 \|gh\mc I(f)\|_{\dot H^2(\Omega)\cap \dot H^4(\Omega)}\lesssim \|g\|_X\|h\|_X,\qquad \|g^2h\mc J(f,g)\|_{\dot H^2(\Omega)\cap \dot H^4(\Omega)}\lesssim \|g\|_X^2\|h\|_X 
 \end{align*}
for all radial $f,g,h\in X_1\cap C^\infty_c(\R^5)$, where $\Omega:=\R^5\setminus\B^5$.

We start with the estimate for $\mc I(f)$. 
By the chain rule,
\begin{align*}
 \langle \cdot\rangle|\nabla \mc I(f)|&\lesssim |\nabla F|+1 \\
 \langle \cdot\rangle|\Delta \mc I(f)|&\lesssim |\Delta F|+|\nabla F|^2+|\nabla F|+1 \\
 \langle \cdot\rangle|\nabla\Delta \mc I(f)|&\lesssim |\nabla\Delta F|
 +|\Delta F||\nabla F|+|\Delta F|+|\nabla F|^3+|\nabla F|^2+|\nabla F|+1 \\
 \langle\cdot\rangle|\Delta^2 \mc I(f)|&\lesssim |\Delta^2 F|
 +|\nabla\Delta F||\nabla F|
 +|\nabla\Delta F|
 +|\Delta F|^2+|\Delta F||\nabla F|^2 \\ &\quad +|\Delta F||\nabla F| 
  +|\Delta F|+|\nabla F|^4+|\nabla F|^3+|\nabla F|^2+|\nabla F|+1,
 \end{align*}
where $F(x)=|x|f(x)$. 
The strategy is to use Lemma \ref{lem:LinfB} to absorb the growing weight in $F$.
We consider $\|gh\Delta \mc I(f)\|_{L^2(\Omega)}$ and estimate
\begin{align*}
\|\langle\cdot\rangle^{-1}gh\Delta F\|_{L^2(\Omega)}&\lesssim \|gh\Delta f\|_{L^2(\Omega)}
+\||\cdot|^{-1}gh\nabla f\|_{L^2(\Omega)}+\||\cdot|^{-2}g h f\|_{L^2(\Omega)} \\
&\lesssim \|g\|_X\|h\|_X\|f\|_X \\
\|\langle\cdot\rangle^{-1}gh |\nabla F|^2\|_{L^2(\Omega)}&\lesssim \||\cdot| gh|\nabla f|^2\|_{L^2(\Omega)}
+\||\cdot|^{-1}ghf^2\|_{L^2(\Omega)} \\
&\lesssim 
\||\cdot|^\frac12 g\|_{L^\infty(\Omega)} \|h\|_{L^\infty(\Omega)}
\||\cdot|^\frac32 \nabla f\|_{L^\infty(\Omega)} 
\||\cdot|^{-1}\nabla f \|_{L^2(\Omega)} \\
&\quad +\||\cdot|^\frac12 g\|_{L^\infty(\Omega)}\||\cdot|^\frac12 h\|_{L^\infty(\Omega)}
\|f\|_{L^\infty(\Omega)}\||\cdot|^{-2}f\|_{L^2(\Omega)} \\
&\lesssim \|g\|_X\|h\|_X\|f\|_X^2 \\
\|\langle\cdot\rangle^{-1}gh \nabla F\|_{L^2(\Omega)}&\lesssim
\|gh\nabla f\|_{L^2(\Omega)}+\||\cdot|^{-1}ghf\|_{L^2(\Omega)} \\
&\lesssim \||\cdot|^\frac12 g\|_{L^\infty(\Omega)}\||\cdot|^\frac12 h\|_{L^\infty(\Omega)}
\left (
\||\cdot|^{-1}\nabla f\|_{L^2(\Omega)}+\||\cdot|^{-2}f\|_{L^2(\Omega)} \right ) \\
&\lesssim \|g\|_X\|h\|_X\|f\|_X
\end{align*}
by Lemma \ref{lem:LinfB} and Hardy's inequality.
This yields $\|\Delta [gh \mc I(f)]\|_{L^2(\Omega)}\lesssim \|g\|_X\|h\|_X$.

Next, we estimate $\|gh\Delta^2 \mc I(f)\|_{L^2(\Omega)}$.
The easy terms are
\begin{align*}
\|\langle\cdot\rangle^{-1}\Delta^2 F\|_{L^2(\Omega)}&\lesssim \sum_{k=0}^4
\||\cdot|^{-4+k}\nabla^k f\|_{L^2(\Omega)}\lesssim \|f\|_X \\
\|\langle\cdot\rangle^{-1}\nabla\Delta F\nabla F\|_{L^2(\Omega)}&\lesssim
\sum_{k=0}^3 \||\cdot|^{-2+k}\nabla^k f\nabla f\|_{L^2(\Omega)}
+\sum_{k=0}^3 \||\cdot|^{-3+k}\nabla^k ff\|_{L^2(\Omega)} \\
&\lesssim \left (\||\cdot|^\frac32 \nabla f\|_{L^\infty(\Omega)}+\|f\|_{L^\infty(\Omega)}\right )\sum_{k=0}^3 \||\cdot|^{-3+k}\nabla^k f\|_{L^2(\Omega)} \\
&\lesssim \|f\|_X^2 \\
\|\langle\cdot\rangle^{-1}\nabla\Delta F\|_{L^2(\Omega)}&\lesssim
\sum_{k=0}^3 \||\cdot|^{-3+k}\nabla^k f\|_{L^2(\Omega)}\lesssim \|f\|_X 
\end{align*}
as well as
\begin{align*}
\|\langle\cdot\rangle^{-1}\Delta F|\nabla F|^2\|_{L^2(\Omega)}&\lesssim 
\||\cdot|^2 \Delta f|\nabla f|^2\|_{L^2(\Omega)}
+\||\cdot| |\nabla f|^3\|_{L^2(\Omega)}+\|f|\nabla f|^2\|_{L^2(\Omega)} \\
&\quad +\||\cdot| \Delta f f^2\|_{L^2(\Omega)}
+\| \nabla f f^2\|_{L^2(\Omega)}+\||\cdot|^{-1}f^3\|_{L^2(\Omega)} \\
&\lesssim \||\cdot|^\frac32 \nabla f\|_{L^\infty(\Omega)}^2
\left (\|\Delta f\|_{L^2(\Omega)}+\||\cdot|^{-1}\nabla f\|_{L^2(\Omega)}
+\||\cdot|^{-2}f\|_{L^2(\Omega)} \right ) \\
&\quad +\||\cdot|^\frac12 f\|_{L^\infty(\Omega)}^2\left (
\|\Delta f\|_{L^2(\Omega)}+\||\cdot|^{-1}\nabla f\|_{L^2(\Omega)}+\||\cdot|^{-2}f\|_{L^2(\Omega)}\right ) \\
&\lesssim \|f\|_X^3
\end{align*}
and
\begin{align*}
\|\langle\cdot\rangle^{-1}\Delta F\nabla F\|_{L^2(\Omega)}&\lesssim
\sum_{k=0}^2\||\cdot|^{-1+k}\nabla^k f\nabla f\|_{L^2(\Omega)}
+\sum_{k=0}^2 \||\cdot|^{-2+k}\nabla^k f f\|_{L^2(\Omega)} \\
&\lesssim \left (\||\cdot|^\frac32 \nabla f\|_{L^\infty(\Omega)}+\|f\|_{L^\infty(\Omega)}\right )\sum_{k=0}^2 \||\cdot|^{-2+k}\nabla^k f\|_{L^2(\Omega)} \\
&\lesssim \|f\|_X^2 \\
\|\langle\cdot\rangle^{-1}|\nabla F|^4\|_{L^2(\Omega)}&\lesssim
\||\cdot|^3|\nabla f|^4\|_{L^2(\Omega)}+\||\cdot|^{-1}f^4\|_{L^2(\Omega)} \\
&\lesssim \||\cdot|^\frac32 \nabla f\|_{L^\infty(\Omega)}^3 \||\cdot|^{-1}\nabla f\|_{L^2(\Omega)}+\||\cdot|^\frac12 f\|_{L^\infty(\Omega)}^3 \||\cdot|^{-2}f\|_{L^2(\Omega)} \\
&\lesssim \|f\|_X^4.
\end{align*}
Analogously, we estimate
\begin{align*}
\|\langle\cdot\rangle^{-1}|\nabla F|^3\|_{L^2(\Omega)}
&\lesssim \||\cdot|^2 |\nabla f|^3\|_{L^2(\Omega)}+\||\cdot|^{-1}f^3\|_{L^2(\Omega)} \\
&\lesssim \||\cdot|^\frac32 \nabla f\|_{L^\infty(\Omega)}^2 \||\cdot|^{-1}\nabla f\|_{L^2(\Omega)}+\||\cdot|^\frac12 f\|_{L^\infty(\Omega)}^2 \||\cdot|^{-2}f\|_{L^2(\Omega)} \\
&\lesssim \|f\|_X^2.
\end{align*}
It remains to control the most delicate term, $|\Delta F|^2$.
For this one we use Hardy and Lemma \ref{lem:LinfB} to obtain
\begin{align*}
\|\langle\cdot\rangle^{-1}gh|\Delta F|^2\|_{L^2(\Omega)}&\lesssim
\||\cdot|gh|\Delta f|^2\|_{L^2(\Omega)}+\||\cdot|^{-1}gh|\nabla f|^2\|_{L^2(\Omega)}
+\||\cdot|^{-3}ghf^2\|_{L^2(\Omega)} \\
&\lesssim \||\cdot|^\frac12 g\|_{L^\infty(\Omega)}
\||\cdot|^\frac12 h\|_{L^\infty(\Omega)} \\
&\quad \times \left (
\||\Delta f|^2\|_{L^2(\Omega)}+\||\cdot|^{-1}|\nabla f|^2\|_{L^2(\Omega)} +\||\cdot|^{-2}f^2\|_{L^2(\Omega)} \right ) \\
&\lesssim \|g\|_X\|h\|_X\|f\|_X^2 .
\end{align*}
The above estimates easily imply $\|\Delta^2[gh \mc I(f)]\|_{L^2(\Omega)}\lesssim \|g\|_X\|h\|_X$.
Putting everything together, we arrive at the desired $\|gh\mc I(f)\|_{\dot H^2(\Omega)\cap\dot H^4(\Omega)}\lesssim \|g\|_X\|h\|_X$.
The bound on $\mc J(f,g)$ is proved in the exact same way.
\end{proof}

The next bound controls the nonlinearity near the center.
Here the issue is to handle powers of $|\cdot|^{-1}$ that arise by differentiation.

\begin{lemma}
\label{lem:MoserB}
Let $\Phi\in C^4(\R)$ and suppose $\Phi'(0)=0$. For $f:\R^5\to\R$ set
\[ \mc M(f)(x)=f(x)^2\Phi(|x|f(x)) \]
Then we have
\[ \|\mc M(f)-\mc M(g)\|_{\dot H^2(\B^5)\cap \dot H^4(\B^5)}\lesssim (\|f\|_X+\|g\|_X)\|f-g\|_X \]
for all $f,g\in X_1 \cap C^\infty_c(\R^5)$.
\end{lemma}

\begin{proof}
As in the proof of Proposition \ref{prop:MoseroB}, we write $\mc I(f)(x)=\Phi(|x|f(x))$ and $\mc I(f)-\mc I(g)=(f-g)\mc J(f,g)$ 
with
\[ \mc J(f,g)(x)=|x|\int_0^1 \Phi'\big (|x|g(x)+t|x|(f(x)-g(x))\big)
\d t \]
and it suffices to show
\[ \|\mc I(f)\|_{\dot H^2(\B^5)\cap \dot H^4(\B^5)}+\|\mc J(f,g)\|_{\dot H^2(\B^5)
\cap \dot H^4(\B^5)}\lesssim 1 \]
for all $f,g\in X_1\cap C^\infty_c(\R^5)$.
We begin with the bound on $\mc I(f)$.
By the chain rule we infer
\begin{align*}
|\mc I(f)|&\lesssim 1 \\
|\nabla \mc I(f)|&\lesssim |\Phi'\circ F||\nabla F| \\
|\Delta \mc I(f)|&\lesssim |\Phi'\circ F||\Delta F|+|\nabla F|^2 \\
|\nabla\Delta \mc I(f)|&\lesssim |\Phi'\circ F||\nabla\Delta F|
+|\Delta F||\nabla F|+|\nabla F|^3 \\
|\Delta^2 \mc I(f)|&\lesssim |\Phi'\circ F||\Delta^2 F|+|\nabla\Delta F||\nabla F|+|\Delta F|^2+
|\Delta F||\nabla F|^2+|\nabla F|^4
\end{align*}
on the ball $\B^5$, where $F(x):=|x|f(x)$ and we have used the fact that $\|F\|_{L^\infty(\B^5)}\lesssim 1$ which follows from Lemma \ref{lem:Linf}.
We consider $\|\Delta \mc I(f)\|_{L^2(\B^5)}$ and estimate
\begin{align*}
\|(\Phi'\circ F) \Delta F\|_{L^2(\B^5)}&\lesssim \|\Delta f\|_{L^2(\B^5)}
+\|\nabla f\|_{L^2(\B^5)}+\||\cdot|^{-1}f\|_{L^2(\B^5)} \\
&\lesssim \|\Delta f\|_{L^2(\R^5)}
+\||\cdot|^{-1}\nabla f\|_{L^2(\R^5)}
+\||\cdot|^{-2}f\|_{L^2(\R^5)}  \\
&\lesssim \|f\|_X \\
\||\nabla F|^2\|_{L^2(\B^5)}&\lesssim 
\||\nabla f|^2\|_{L^2(\B^5)}+\|f^2\|_{L^2(\B^5)} \\
&\lesssim \|\nabla f\|_{L^\infty(\R^5)}\||\cdot|^{-1}\nabla f\|_{L^2(\R^5)}
+\|f\|_{L^\infty(\R^5)}\||\cdot|^{-2}f\|_{L^2(\R^5)} \\
&\lesssim \|f\|_X^2 
\end{align*}
by Lemma \ref{lem:Linf} and Hardy's inequality.
This yields the desired $\|\Delta \mc I(f)\|_{L^2(\B^5)}\lesssim 1$. 

Next, we estimate $\|\Delta^2 \mc I(f)\|_{L^2(\B^5)}$. The most delicate term is
$|\Phi'\circ F||\Delta^2 F|$ where we absorb one singular factor $|\cdot|^{-1}$ by exploiting the assumption $\Phi'(0)=0$.
More precisely,
\begin{align*}
\|(\Phi'\circ F) \Delta^2 F\|_{L^2(\B^5)}
&\lesssim \||\cdot|^{-1}\Phi'\circ F\|_{L^\infty(\B^5)}\sum_{k=0}^4 \|\cdot|^{-2+k}\nabla^k f\|_{L^2(\B^5)} \\
&\lesssim \|f\|_X^2
\end{align*}
by Hardy's inequality, the bound
\[ |(\Phi'\circ F)(x)|=|\Phi'(|x|f(x))|\lesssim |x||f(x)|, \]
which follows from $\Phi'(0)=0$, and Lemma \ref{lem:Linf}.
Furthermore, by Gagliardo-Nirenberg (see e.g.~\cite{Tay11}, p.~8, Proposition 3.1),
\begin{align*}
\||\Delta F|^2\|_{L^2(\B^5)}&\lesssim \|(\chi \Delta F)^2\|_{L^2(\R^5)} 
\lesssim \|\nabla (\chi\nabla F)\nabla(\chi\nabla F)\|_{L^2(\B^5)}
+\|(\nabla\chi\nabla F)^2\|_{L^2(\R^5)} \\
&\lesssim \|\nabla(\chi\nabla F)\|_{L^4(\R^5)}^2
+\|(\nabla\chi\nabla F)^2\|_{L^2(\R^5)} \\
&\lesssim \|\chi\nabla F\|_{L^\infty(\R^5)}\|\Delta(\chi \nabla F)\|_{L^2(\R^5)} 
+\|\nabla\chi \nabla F\|_{L^\infty(\R^5)}^2 \\
&\lesssim \|f\|_X^2,
\end{align*}
where $\chi: \R^5\to [0,1]$ is a smooth cut-off satisfying $\chi(x)=1$ for $|x|\leq 1$ and $\chi(x)=0$ for $|x|\geq 2$.
The remaining terms are readily estimated as
\begin{align*}
\|\nabla\Delta F\nabla F\|_{L^2(\B^5)}&\lesssim \|\nabla F\|_{L^\infty(\B^5)}
\|\nabla \Delta F\|_{L^2(\B^5)}\lesssim \|f\|_X^2 \\
\|\Delta F|\nabla F|^2\|_{L^2(\B^5)}&\lesssim \|\nabla F\|_{L^\infty(\B^5)}^2
\|\Delta F\|_{L^2(\B^5)}\lesssim \|f\|_X^3
\end{align*}
This shows $\|\Delta^2\mc I(f)\|_{L^2(\B^5)}\lesssim 1$.
The proof of the bound on $\mc J(f,g)$ is identical.
\end{proof}

In fact, we need a slightly more general form of Lemma \ref{lem:MoserB}.

\begin{corollary}
\label{cor:MoserB}
Let $\Phi\in C^4(\R^2)$ and suppose $\partial_1\Phi(0,0)=\partial_2\Phi(0,0)=0$. Set
\[ \mc M(f)(x)=f(x)^2\Phi(|x|f(x),|x|\varphi_0(x)) \]
where $\varphi_0\in C^\infty_c(\R^5)$ is a fixed function. 
Then we have
\[ \|\mc M(f)-\mc M(g)\|_{\dot H^2(\B^5)\cap \dot H^4(\B^5)}\lesssim (\|f\|_X+\|g\|_X)\|f-g\|_X \]
for all $f\in X_1 \cap C^\infty_c(\R^5)$.
\end{corollary}

\begin{proof}
This is a straightforward generalization of Lemma \ref{lem:MoserB}.
\end{proof}

We are now in a position to prove that the nonlinearity $\mc N$ is locally Lipschitz on $X$.

\begin{lemma}
\label{lem:N}
We have the bound
\[ \|\mc N(f)-\mc N(g)\|\lesssim (\|f\|_X+\|g\|_X)\|f-g\|_X \]
for all $f,g \in X_1$.
\end{lemma}

\begin{proof}
By a density argument it suffices to consider $f,g\in X_1\cap C^\infty_c(\R^5)$.
Recall that 
$\mc N(f)(x)=N(f(x),|x|)$ with
\[ N(u,y)=-\frac{1}{y^3}\left [
\sin(2f_0(y)+2yu)-\sin(2f_0(y))-2y\cos(2f_0(y))u \right ].
\]
Note that
\begin{align*}
\partial_1 N(u,y)&=-\frac{2}{y^2}\left [
\cos(2f_0(y)+2yu)-\cos(2f_0(y)) \right ] \\
\partial_1^2 N(u,y)&=\frac{4}{y}\sin(2f_0(y)+2yu) \\
\partial_1^3 N(u,y)&=8\cos(2f_0(y)+2yu) .
\end{align*}
Evidently, $N(0,y)=\partial_1 N(0,y)=0$ and the fundamental theorem of calculus yields
\begin{align*}
 N(u,y)&=\int_0^1 \partial_{t_1}N(t_1u,y)\d t_1 = u\int_0^1 \partial_1 N(t_1 u,y)\d t_1  
=u^2 \int_0^1 t_1 \int_0^1 \partial_1^2 N(t_2t_1 u,y)\d t_2 \d t_1 \\
&=\frac{4u^2}{y}\int_0^1 t_1 \int_0^1 \sin(2f_0(y)+2t_2t_1yu)\d t_2 \d t_1 \\
&=8u^2 \frac{f_0(y)}{y}\int_0^1 t_1\int_0^1 \int_0^1 \cos(2t_3f_0(y)+2t_3t_2t_1yu)\d t_3 \d t_2 \d t_1 \\
&\quad +8u^3\int_0^1 t_1^2\int_0^1 t_2\int_0^1 \cos(2t_3f_0(y)+2t_3t_2t_1yu)\d t_3 \d t_2 \d t_1.
 \end{align*}
 We define $\Phi_1, \tilde \Phi_1: \R^2\to\R$ and $\Phi_2: \R\times \R^5\to \R$ by
 \begin{align*}
  \Phi_1(v,v_0)&:=8\int_0^1 t_1\int_0^1 \int_0^1 \cos(2t_3v_0+2t_3t_2t_1v)\d t_3 \d t_2 \d t_1 \\
   \tilde \Phi_1(v,v_0)&:=8\int_0^1 t_1^2\int_0^1 t_2\int_0^1 \cos(2t_3v_0+2t_3t_2t_1v)\d t_3 \d t_2 \d t_1 \\
  \Phi_2(v,x)&:=\frac{4}{|x|}\int_0^1 t_1 \int_0^1 \sin(2f_0(|x|)+2t_2t_1v)\d t_2 \d t_1
  \end{align*}
  and set
  \begin{align*} \mc M_1(f)(x)&:=f(x)^2 \Phi_1(|x|f(x),|x|\varphi_0(x)) \\
  \tilde{\mc M}_1(f)(x)&:=f(x)^2 \tilde \Phi_1(|x|f(x),|x|\varphi_0(x)) \\
  \mc M_2(f)(x)&:=f(x)^2 \Phi_2(|x|f(x),x)  
  \end{align*}
  where $\varphi_0(x):=\chi(x)f_0(|x|)/|x|$ with $\chi:\R^5\to[0,1]$ the usual smooth cut-off satisfying $\chi(x)=1$ for $|x|\leq 1$ and $\chi(x)=0$ for $|x|\geq 2$.
  By \cite{BieDon16}, $f_0$ is odd and thus, $\varphi_0\in C^\infty_c(\R^5)$.
  This yields the representations
  \[ \mc N(f)(x)=\varphi_0(x)\mc M_1 f(x)+f(x)\tilde{\mc M}_1 (f)(x) \]
  for $|x|\leq 1$ and
  \[ \mc N(f)(x)=\mc M_2(f)(x) \]
  for $|x|\geq \frac12$.
Evidently, we have 
\[ \partial_1 \Phi_1(0,0)=\partial_2\Phi_1(0,0)=\partial_1 \tilde \Phi_1(0,0)=\partial_2\tilde \Phi_1(0,0)=0 \] and
 \[ |\partial^\alpha \Phi_2(v,x)|\leq C_\alpha \langle x\rangle^{-1} \]
 for all $(v,x)\in \R\times \R^5\setminus\B^5$ and all multi-indices $\alpha\in \N_0^6$.
 As a consequence, Corollary \ref{cor:MoserB} and Proposition \ref{prop:MoseroB} apply to $\mc M_1$, $\tilde{\mc M}_1$, and $\mc M_2$, respectively, and Corollary \ref{cor:alg} yields the claim.
\end{proof}

\subsection{The initial data operator}
Now we consider the initial data operator
\[ \mc U(h,T)(x):=f_0(\sqrt T |x|)/|x|-f_0(|x|)/|x|+h(\sqrt T |x|)/|x| .
 \]
Recall that
\[ \tilde Y=\{h\in C^\infty_c([0,\infty)): h^{(2k)}(0)=0 \mbox{ for all }k\in \N_0\} \]
and $Y$ was defined as the completion of $\tilde Y$ with respect to the norm
\[ \|h\|_Y=\||\cdot|^{-1}h(|\cdot|)\|_X. \]
First, we need to make sure that $\mc U(h,T)$ has values in $X$. The following more general result will be helpful in this respect.

\begin{lemma}
\label{lem:finX}
Let $f\in C^\infty(\R^5)$ and assume the bounds
\[ |\nabla^k f(x)|\lesssim \langle x\rangle^{-1-k} \]
for all $x\in \R^5$ and $k\in \{0,1,2,3,4\}$. Then $f\in X$.
\end{lemma}

\begin{proof}
Let $\chi: \R^5\to [0,1]$ be the usual smooth cut-off satisfying $\chi(x)=1$ for $|x|\leq 1$ and $\chi(x)=0$ for $|x|\geq 2$. For $n\in \N$ we set $\chi_n(x):=\chi(x/n)$.
Then $\chi_n f\in C^\infty_c(\R^5)$ for any $n\in \N$ and $\chi_nf -f=0$ on $\B^5_n$.
Thus, thanks to the decay $|f(x)|\lesssim \langle x \rangle^{-1}$, 
\begin{align*} \|\chi_n f-f\|_{L^\infty(\R^5)}&\leq \|\chi_n f-f\|_{L^\infty(\B^5_n)}
+\|\chi_n f-f\|_{L^\infty(\R^5\setminus\B^5_n)}\lesssim \|f\|_{L^\infty(\R^5\setminus \B^5_n)} \\
&\lesssim n^{-1} 
\end{align*}
and we see that $\chi_n f\to f$ in $L^\infty(\R^5)$.
Now let $m\leq n$ and note that
$\chi_n f-\chi_m f=0$ on the ball $\B^5_m$. 
Furthermore,
\begin{align*}
 \|\chi_n f\|_{\dot H^2(\R^5\setminus \B^5_m)}&\simeq
 \|\Delta\chi_n f\|_{L^2(\R^5\setminus \B^5_m)}
 +\|\nabla\chi_n \nabla f\|_{L^2(\R^5\setminus \B^5_m)}
 +\|\chi_n \Delta f\|_{L^2(\R^5\setminus \B^5_m)} \\
 &\lesssim \||\cdot|^{-1}\Delta\chi_n\|_{L^2(\R^5\setminus \B^5_m)}
 +\||\cdot|^{-2}\nabla\chi_n\|_{L^2(\R^5\setminus \B^5_m)}
 +\||\cdot|^{-3}\chi_n\|_{L^2(\R^5\setminus \B^5_m)} \\
 &\lesssim m^{-\frac12}
 \end{align*}
and similarly for $\|\chi_n f\|_{\dot H^4(\R^5\setminus \B^5_m)}$.
In summary, we find
\[ \|\chi_mf -\chi_n f\|_X\lesssim m^{-\frac12}+n^{-\frac12} \]
for all $n,m\in \N$ and thus,
$(\chi_n f)_{n\in\N}$ is Cauchy in $X$. This shows $f\in X$, cf.~Lemma \ref{lem:embed}.
\end{proof}

\begin{corollary}
\label{cor:f0inX}
The function $x\mapsto f_0(|x|)/|x|: \R^5\to\R$ belongs to $X$.
\end{corollary}

\begin{proof}
By \cite{BieDon16}, $f_0: [0,\infty)\to\R$ is smooth, odd, and satisfies the bounds $|f_0^{(k)}(y)|\lesssim y^{-k}$ for all $y\geq 1$ and $k\in \{0,1,2,3,4\}$. Consequently, the function $x\mapsto f_0(|x|)/|x|: \R^5\to\R$ verifies the hypotheses of Lemma \ref{lem:finX}.
\end{proof}

\begin{lemma}
\label{lem:U}
The map $\mc U: Y\times [\frac12,\frac32]\to X$ is well-defined and continuous.
Furthermore, we have the bound
\[ \|\mc U(h,T)\|_X\lesssim \|h\|_Y+|T-1| \]
for all $(h,T)\in Y\times [\frac12,\frac32]$.
\end{lemma}

\begin{proof}
Corollary \ref{cor:f0inX} and the very definition of $Y$ show that $\mc U$ has values in $X$. Hence, $\mc U$ is well-defined.
For brevity we set $I:=[\frac12,\frac32]$, $\varphi_0(x):=f_0(|x|)/|x|$, and $H(x):=h(|x|)/|x|$.
Since $\|h\|_Y=\|H\|_X$, 
we have to show that the maps $T \mapsto \varphi_0(\sqrt T(\cdot)): I\to X$ and $(H,T)\mapsto H(\sqrt T(\cdot)): X\times I\to X$ are continuous.
The fundamental theorem of calculus yields
\begin{align*} \varphi_0(\sqrt{T_1} x)-\varphi_0(\sqrt{T_2} x)&=\int_0^1 
\partial_t \varphi_0\left (\sqrt{T_2} x+t\left (\sqrt{T_1}-\sqrt{T_2}\right )x\right )\d t  \\
&=\left (\sqrt{T_1}-\sqrt{T_2}\right )\int_0^1 x\nabla \varphi_0\left (\sqrt{T_2}x+t\left (\sqrt{T_1}-\sqrt{T_2}\right )x\right )\d t \\
&=:\left (\sqrt{T_1}-\sqrt{T_2}\right )\psi_{T_1,T_2}(x).
\end{align*}
Obviously, $\psi_{T_1,T_2}\in C^\infty(\R^5)$ and $|\nabla^k \psi_{T_1,T_2}(x)|\lesssim \langle x\rangle^{-1-k}$ for all $x\in \R^5$ and $k\in \{0,1,2,3,4\}$.
Consequently, $\psi_{T_1,T_2}\in X$ by Lemma \ref{lem:finX} and we infer
\[ \left \|\varphi_0(\sqrt{T_1}(\cdot))-\varphi_0(\sqrt{T_2}(\cdot))\right \|_X\leq \left |\sqrt{T_1}-\sqrt{T_2}\right |\|\psi_{T_1,T_2}\|_X\lesssim \left |\sqrt{T_1}-\sqrt {T_2}\right | \]
for all $T_1,T_2\in I$.

Now let $\epsilon>0$ and choose $H_1, H_2 \in X$ such that $\|H_1-H_2\|_X\leq \epsilon/100$. Furthermore, choose $\tilde H_1 \in C^\infty_c(\R^5)$ such that $\|H_1-\tilde H_1\|_X\leq \epsilon/100$.
Then we have
\begin{align*}
\left \|H_1(\sqrt{T_1}(\cdot))-H_2(\sqrt{T_2}(\cdot))\right \|_X&\leq
\left \|H_1(\sqrt{T_1}(\cdot))-H_1(\sqrt{T_2}(\cdot))\right \|_X \\
&\quad +\left \|H_1(\sqrt{T_2}(\cdot))-H_2(\sqrt{T_2}(\cdot))\right \|_X \\
&\leq \left  \|H_1(\sqrt{T_1}(\cdot))-H_1(\sqrt{T_2}(\cdot))\right \|_X +\tfrac{\epsilon}{4}
\end{align*}
and
\begin{align*}
 \left \|H_1(\sqrt{T_1}(\cdot))-H_1(\sqrt{T_2}(\cdot))\right \|_X&\leq
\left \|H_1(\sqrt{T_1}(\cdot))-\tilde H_1(\sqrt{T_1}(\cdot))\right \|_X \\
&\quad +\left \|\tilde H_1(\sqrt{T_1}(\cdot))-\tilde H_1(\sqrt{T_2}(\cdot))\right \|_X\\
&\quad +\left \|\tilde H_1(\sqrt{T_2}(\cdot))-H_1(\sqrt{T_2}(\cdot))\right \|_X \\
&\leq  \tfrac{\epsilon}{2}+C_\epsilon \left |\sqrt{T_1}-\sqrt{T_2}\right |
\end{align*}
for all $T_1,T_2\in I$ again by the fundamental theorem of calculus.
Consequently, we may choose $|T_1-T_2|$ so small that 
\[ \left \|H_1(\sqrt{T_1}(\cdot))-H_2(\sqrt{T_2}(\cdot))\right \|_X<\epsilon. \]
This proves the continuity of $\mc U$.
Finally, from the above it is obvious that
\[ \|\mc U(h,T)\|_X\lesssim \|h\|_Y+|\sqrt T-1|\lesssim \|h\|_Y+|T-1|. \]
\end{proof}

\subsection{Global existence for the modified equation}
Now we turn to the solution of Eq.~\eqref{eq:fullop}.
As an intermediate step we consider the Cauchy problem
\begin{equation}
\label{eq:fullopint}
\left \{ \begin{array}{l}
\partial_s \phi(s)=\mc L\phi(s)+\mc N(\phi(s)) \\
\phi(0)=f
\end{array} \right. 
\end{equation}
for given small $f\in X$.
We employ Duhamel's principle to obtain the weak formulation
\begin{equation}
\label{eq:fullintweak}
 \phi(s)=e^{s\mc L}f+\int_0^s e^{(s-s')\mc L}\mc N(\phi(s'))\d s'. 
 \end{equation}
As a matter of fact, this equation does not have global solutions for arbitrary $f$ due to the unstable subspace of the semigroup $e^{s\mc L}$.
Thus, we modify Eq.~\eqref{eq:fullintweak} by adding a correction term that stabilizes the evolution.
In order to obtain this term, we formally project the evolution to the unstable subspace.
That is to say, we define the projection operator $\mc P: H\to H$ by 
\[ \mc Pf:=(f|\psi_1)_{L^2_\sigma(\R^5)}\psi_1. \]
Note that by \cite{BieDon16}, $\psi_1(x)=f_0'(|x|)/\|f_0'(|\cdot|)\|_{L^2_\sigma(\R^5)}$ satisfies the assumptions of Lemma \ref{lem:finX} and thus,
$\mc P$ has values in $X$.
Furthermore, by Lemma \ref{lem:Linf}, $\mc P|_X$ is a bounded projection on $X$.
Applying $\mc P$ to Eq.~\eqref{eq:fullintweak} and using the fact that $\mc Pe^{s\mc L}f=e^s \mc Pf$, we obtain (at least formally)
\[ \mc P\phi(s)=e^s \mc Pf+e^s \int_0^s e^{-s'}\mc P\mc N(\phi(s'))\d s'. \]
This suggests to subtract the term $e^s \mc C(\Phi, f)$, where
\[ \mc C(\phi, f):=\mc Pf+\int_0^\infty e^{-s'}\mc P\mc N(\phi(s'))\d s'. \]
In order to put this on a sound functional analytic footing, we introduce the Banach space
\[ \mc X:=\left \{\phi \in C([0,\infty),X): \|\phi\|_{\mc X}<\infty \right \} \]
with the norm
\[ \|\phi\|_{\mc X}:=\sup_{s>0} e^{\omega_0 s}\|\phi(s)\|_X, \]
where $\omega_0>0$ is the constant from Theorem \ref{thm:lin}.
Furthermore, for $\delta>0$, we set
\[ \mc X_\delta:=\{\phi\in \mc X: \|\phi\|_{\mc X}\leq \delta\}. \]
Now we define $\mc K: \mc X\times X\to \mc X$ by
\[ \mc K(\phi, f)(s):=e^{s\mc L}f+\int_0^s e^{(s-s')\mc L}\mc N(\phi(s'))\d s'
-e^s \mc C(\phi,f) \]
and show that $\mc K(\cdot,f)$ is a contraction on $\mc X_\delta$, provided $f\in X$ is sufficiently small.

\begin{lemma}
\label{lem:self}
There exists a constant $c>0$ such that $\mc K(\phi,f)\in \mc X_\delta$ for all $\phi\in \mc X_\delta$ and all $f\in X_{\delta/c}$, provided $\delta>0$ is sufficiently small.
\end{lemma}

\begin{proof}
By definition, we have
\[ \mc P\mc K(\phi,f)(s)=-e^s\int_s^\infty e^{-s'}\mc P\mc N(\phi(s'))\d s' \]
and thus, by Lemma \ref{lem:N},
\begin{align*} \|\mc P\mc K(\phi,f)(s)\|_X&\lesssim
e^s \int_s^\infty e^{-s'}\|\phi(s')\|_X^2 \d s'\lesssim e^{-2\omega_0 s}\|\phi\|_{\mc X}^2 \lesssim \delta^2 e^{-2\omega_0 s}. 
\end{align*}
Similarly, 
\[ (1-\mc P)\mc K(\phi,f)(s)=e^{s\mc L}(1-\mc P)f+\int_0^s e^{(s-s')\mc L}
(1-\mc P)\mc N(\phi(s'))\d s' \]
and thus,
\begin{align*}
\|(1-\mc P)\mc K(\phi,f)(s)\|_X&\lesssim e^{-\omega_0 s}\|f\|_X+\int_0^s e^{-\omega_0 (s-s')}\|\phi(s')\|_X^2 \d s' \\
&\lesssim \tfrac{\delta}{c}e^{-\omega_0 s}+\|\phi\|_{\mc X}^2 e^{-\omega_0 s}\int_0^s e^{-\omega_0 s'}\d s' \\
&\lesssim \tfrac{\delta}{c}e^{-\omega_0 s}+\delta^2 e^{-\omega_0 s}
\end{align*}
by Theorem \ref{thm:lin} and Lemma \ref{lem:N}.
In summary, this yields $\|\mc K(\phi,f)\|_{\mc X}\lesssim \frac{\delta}{c}+\delta^2$ and by choosing $c>0$ large enough and $\delta>0$ small enough, we obtain
$\|\mc K(\phi,f)\|_{\mc X}\leq \delta$.
\end{proof}

\begin{lemma}
\label{lem:contr}
Let $\delta>0$ be sufficiently small. Then we have the bound
\[ \|\mc K(\phi,f)-\mc K(\psi,f)\|_{\mc X}\leq \tfrac12 \|\phi-\psi\|_{\mc X} \]
for all $\phi,\psi\in \mc X_\delta$ and all $f\in X$.
\end{lemma}

\begin{proof}
We have
\[ \mc P\mc K(\phi,f)(s)-\mc P\mc K(\psi,f)(s)=-e^s\int_s^\infty e^{-s'}
\mc P\big [ \mc N(\phi(s'))-\mc N(\psi(s')) \big ]\d s' \]
and Lemma \ref{lem:N} yields
\[ \|\mc P\mc K(\phi,f)(s)-\mc P\mc K(\psi,f)(s)\|_X\lesssim \delta e^{-2\omega_0 s}\|\phi-\psi\|_{\mc X}, \]
cf.~the proof of Lemma \ref{lem:self}.
Similarly,
\[ (1-\mc P)\mc K(\phi,f)(s)-(1-\mc P)\mc K(\psi,f)(s)=
\int_0^s e^{(s-s')\mc L}(1-\mc P)\big [ \mc N(\phi(s'))-\mc N(\psi(s')) \big ]\d s'  \]
and thus,
\[ \|(1-\mc P)\mc K(\phi,f)(s)-(1-\mc P)\mc K(\psi,f)(s)\|_X\lesssim \delta e^{-\omega_0 s}\|\phi-\psi\|_{\mc X}. \]
In summary, we infer
\[ \|\mc K(\phi,f)-\mc K(\psi,f)\|_{\mc X}\lesssim \delta \|\phi-\psi\|_{\mc X} \]
and by choosing $\delta>0$ sufficiently small, we arrive at the claim.
\end{proof}

Based on the above, it is now easy to construct a global solution to the \emph{modified equation}

\begin{equation}
\label{eq:mod} \phi(s)=e^{s\mc L}\mc U(h,T)+\int_0^s e^{(s-s')\mc L}\mc N(\phi(s'))\d s'-e^s\mc C(\phi,\mc U(h,T)). 
\end{equation}

\begin{corollary}
\label{cor:ex}
Let $M>0$ be sufficiently large and $\delta>0$ sufficiently small. Then, for every $h\in Y$ and every $T>0$ satisfying
\[ \|h\|_Y+|T-1|\leq \tfrac{\delta}{M}, \]
there exists a unique $\phi_{h,T}\in \mc X_\delta$ such that 
\[ \phi_{h,T}=\mc K(\phi_{h,T},\mc U(h,T)). \]
In particular, $\phi_{h,T}$ is a solution to Eq.~\eqref{eq:mod}.
Furthermore, the solution map $(h,T)\mapsto \phi_{h,T}$ is continuous.
\end{corollary}

\begin{proof}
By Lemma \ref{lem:U}, we can achieve
\[ \|\mc U(h,T)\|_X\leq \tfrac{\delta}{c} \]
for any given $c>0$ by choosing $M$ sufficiently large. 
Thus, the existence and uniqueness of $\phi_{h,T}$ is a consequence of Lemmas \ref{lem:self}, \ref{lem:contr}, and the contraction mapping principle.

For the continuity of the solution map we note that
\begin{align*}
\|\phi_{h_1,T_1}-\phi_{h_2,T_2}\|_{\mc X}&=\|\mc K(\phi_{h_1,T_1},\mc U(h_1,T_1))
-\mc K(\phi_{h_2,T_2},\mc U(h_2,T_2))\|_{\mc X} \\
&\leq \|\mc K(\phi_{h_1,T_1},\mc U(h_1,T_1))
-\mc K(\phi_{h_2,T_2},\mc U(h_1,T_1))\|_{\mc X} \\
&\quad +\|\mc K(\phi_{h_2,T_2},\mc U(h_1,T_1))
-\mc K(\phi_{h_2,T_2},\mc U(h_2,T_2))\|_{\mc X} \\
&\leq \tfrac12 \|\phi_{h_1,T_1}-\phi_{h_2,T_2}\|_{\mc X}+C\|\mc U(h_1,T_1)-\mc U(h_2,T_2)\|_X
\end{align*}
by Lemma \ref{lem:contr} since
\begin{align*}
\|\mc K&(\phi_{h_2,T_2},\mc U(h_1,T_1))(s)
-\mc K(\phi_{h_2,T_2},\mc U(h_2,T_2))(s)\|_X \\
&=\|e^{s\mc L}(1-\mc P)[\mc U(h_1,T_1)-\mc U(h_2,T_2)]\|_X \\
&\lesssim e^{-\omega_0 s}\|\mc U(h_1,T_1)-\mc U(h_2,T_2)\|_X
\end{align*}
by Theorem \ref{thm:lin}.
Consequently, Lemma \ref{lem:U} finishes the proof.
\end{proof}

Corollary \ref{cor:ex} provides us with a solution to the modified equation
\eqref{eq:mod}.
Thus, in order to obtain a (mild) solution to Eq.~\eqref{eq:fullop}, we have to get rid of the correction term $\mc C(\phi,\mc U(h,T))$.
So far, $h$ and $T$ can be chosen freely, subject to the smallness conditions in Corollary \ref{cor:ex}.
In the last step of the construction we now show that for any small $h\in Y$ there exists in fact a $T_h>0$ such that $\mc C(\phi_{h,T_h}, \mc U(h,T_h))=0$.

\begin{lemma}
\label{lem:corr}
Let $M>0$ be sufficiently large and $\delta>0$ sufficiently small. 
Then, for every $h\in Y$ satisfying $\|h\|_Y\leq \frac{\delta}{M^2}$, there exists a $T_h \in [1-\frac{\delta}{M}, 1+\frac{\delta}{M}]$ such that
\[ \mc C(\phi_{h,T_h},\mc U(h,T_h))=0. \]
\end{lemma}

\begin{proof}
For brevity we set $I_{M,\delta}=[1-\frac{\delta}{M}, 1+\frac{\delta}{M}]$.
The map $\mc C$ has values in $\langle \psi_1\rangle$ and thus, it suffices to consider the real-valued function $F_h: I_{M,\delta}\to \R$ given by
\[ F_h(T):=\big (\mc C(\phi_{h,T},\mc U(h,T)) \big | \psi_1 \big )_{L^2_\sigma(\R^5)}. \]
By Corollary \ref{cor:ex}, $F_h$ is continuous.
Furthermore, by noting that
\[ \left .\partial_T \frac{f_0(\sqrt T|x|)}{|x|}\right|_{T=1}=\tfrac12 f_0'(|x|)=\tfrac 12 \|f_0(|\cdot|)\|_{L^2_\sigma(\R^5)}\psi_1(x) \]
we obtain by a Taylor expansion the representation
\[ \mc U(h,T)(x)=\tfrac12(T-1)\|f_0(|\cdot|)\|_{L^2_\sigma(\R^5)}\psi_1(x)+(T-1)^2 f_T(x) + h(\sqrt T|x|)/|x| \]
where $T\mapsto f_T: I_{M,\delta}\to X$ is continuous and $\|f_T\|_X\lesssim 1$ for all $T\in I_{M,\delta}$.
This yields
\begin{align*}
 F_h(T)&=\big (\mc C(\phi_{h,T},\mc U(h,T)) \big | \psi_1 \big )_{L^2_\sigma(\R^5)} \\
&=\big (\mc P\mc U(h,T)\big |\psi_1\big )_{L^2_\sigma(\R^5)}
+\int_0^\infty e^{-s'}\big (\mc P\mc N(\phi_{h,T}(s'))\big | \psi_1\big )_{L^2_\sigma(\R^5)}\d s' \\
&=\tfrac12 (T-1)\|f_0(|\cdot|)\|_{L^2_\sigma(\R^5)}+O(\tfrac{\delta}{M^2} T^0)+O(\delta^2 T^0).
\end{align*}
Consequently, by setting $\tilde F_h(T)=2\|f_0(|\cdot|)\|_{L^2_\sigma(\R^5)}^{-1}F_h(T)$, we infer
\[ \tilde F_h(T)=T-1+O(\tfrac{\delta}{M^2}T^0)+O(\delta^2T^0) \]
and $\tilde F_h(T)=0$ is equivalent to $T-1=G(T)$ for a continuous function $G: I_{M,\delta}\to \R$ that satisfies
\[ |G(T)|\leq C\tfrac{\delta}{M^2}+C\delta^2 \]
for all $T\in I_{M,\delta}$.
By choosing $M>0$ sufficiently large and $\delta>0$ sufficiently small, we can achieve $|G(T)|\leq \frac{\delta}{M}$ for all $T\in I_{M,\delta}$ and thus, $1+G$ is a continuous self-map of the interval $I_{M,\delta}$ which necessarily has a fixed point $T_h\in I_{M,\delta}$.
\end{proof}

\subsection{Proof of Theorem \ref{thm:main}}
Without loss of generality we set $T_0=1$.
Lemma \ref{lem:corr} yields a strong solution $w_h(y,s)=\phi_{h,T_h}(s)(ye_1)$ of Eq.~\eqref{eq:full} with $T=T_h$ and
\[ \|w_h(|\cdot|,s)\|_X\leq \delta e^{-\omega_0 s}. \]
By construction, see Section \ref{sec:pre}, 
\begin{align*}
 u_h(r,t):=f_0\left (\frac{r}{\sqrt{T_h-t}}\right )+\frac{r}{\sqrt{T_h-t}}w_h\left (\frac{r}{\sqrt{T_h-t}}, -\log(T_h-t)+\log T_h\right )
\end{align*}
is a solution of Eq.~\eqref{eq:hmhf} with initial data $u_h(r,0)=u_1^*(0,r)+h(r)$.
By scaling, we infer
\begin{align*}
\|u_h(\cdot,t)-u_{T_h}^*(\cdot,t)\|_Y
&=(T_h-t)^{-\frac12}\left \|w_h\left (\frac{|\cdot|}{\sqrt{T_h-t}}, -\log(T_h-t)+\log T_h\right )\right \|_X \\
&\lesssim (T_h-t)^{-\frac54}\left \|w_h(|\cdot|,-\log(T_h-t)+\log T_h)\right \|_X \\
&\lesssim \delta (T_h-t)^{-\frac54+\omega_0}
\end{align*}
for all $t\in [0,T_h)$.
Furthermore, from Corollary \ref{cor:f0inX} it follows that $f_0$ belongs to $Y$ and the blowup speed of $u_{T_h}^*$ in $Y$ is
\begin{align*}
 \|u_{T_h}^*(\cdot,t)\|_Y&=\left \||\cdot|^{-1}f_0\left (\frac{|\cdot|}{\sqrt{T_h-t}}\right ) \right \|_X \\
&=(T_h-t)^{-\frac14}\||\cdot|^{-1}f_0(|\cdot|)\|_{\dot H^2(\R^5)} 
+(T_h-t)^{-\frac54}\||\cdot|^{-1}f_0(|\cdot|)\|_{\dot H^4(\R^5)} \\
&\simeq (T_h-t)^{-\frac54}. 
 \end{align*}
Consequently, the statement of Theorem \ref{thm:main} follows by choosing $M$ sufficiently large.

\bibliography{hmhf}
\bibliographystyle{plain}

\end{document}